\documentclass[11 pt]{amsart}
\usepackage{thmtools, thm-restate}
\usepackage{amsmath, amsthm,amssymb,bbm,enumerate,mathrsfs,mathtools}
\usepackage{a4wide}
\usepackage{tikz,xcolor,verbatim}
\usetikzlibrary{calc,shapes}
\usepackage{hyperref}
\usepackage{makecell}
\usepackage{mathrsfs}
\usepackage[numbers,sort&compress]{natbib}
\usepackage{cancel}
\usepackage{caption}
\usetikzlibrary{shapes.geometric}

\tikzset{
triangle/.style={
    draw,
    shape border rotate=0,
    regular polygon,
    regular polygon sides=3,
    fill=white,
    node distance=2cm,
    minimum height=4em
}
}
\usepackage{amssymb}
\usepackage{dirtytalk}
\usepackage{algorithmicx}
\usepackage{algpseudocode}
\usepackage{algorithm}
\usepackage{a4wide}
\algblock{Input}{EndInput}
\algnotext{EndInput}
\algblock{Output}{EndOutput}
\algnotext{EndOutput}
\algblock{Initialization}{EndInitialization}
\algnotext{EndInitialization}
\algblock{Begin}{EndBegin}
\algnotext{EndBegin}

\usepackage[yyyymmdd,hhmmss]{datetime}
\usepackage{background}
\SetBgContents{ver. \currenttime \qquad \today }
\SetBgScale{1}
\SetBgAngle{0}
\SetBgOpacity{1}
\SetBgPosition{current page.south west} 
\SetBgHshift{3cm}
\SetBgVshift{0.5cm} 

\tikzset{
bigblue/.style={circle, draw=blue!80,fill=blue!40,thick, inner sep=1.5pt, minimum size=5mm},
bigred/.style={circle, draw=red!80,fill=red!40,thick, inner sep=1.5pt, minimum size=5mm},
bigblack/.style={circle, draw=black!100,fill=black!40,thick, inner sep=1.5pt, minimum size=5mm},
bluevertex/.style={circle, draw=blue!100,fill=blue!100,thick, inner sep=0pt, minimum size=2mm},
redvertex/.style={circle, draw=red!100,fill=red!100,thick, inner sep=0pt, minimum size=2mm},
blackvertex/.style={circle, draw=black!100,fill=black!100,thick, inner sep=0pt, minimum size=1.5mm},  
whitevertex/.style={circle, draw=black!100,fill=white!100,thick, inner sep=0pt, minimum size=2mm},  
smallblack/.style={circle, draw=black!100,fill=black!100,thick, inner sep=0pt, minimum size=1mm},
smallwhite/.style={circle, draw=white!100,fill=white!100,thick, inner sep=0pt, minimum size=1mm},
redvertexv2/.style={circle, draw=black!100,fill=red!100,thick, inner sep=0pt, minimum size=3.35mm},
bluevertexv2/.style={circle, draw=black!100,fill=blue!100,thick, inner sep=0pt, minimum size=3.35mm},
greenvertexv2/.style={circle, draw=black!100,fill=green!100,thick, inner sep=0pt, minimum size=3.35mm},
blackvertexv2/.style={circle, draw=black!100,fill=black!100,thick, inner sep=0pt, minimum size= 3.35mm},
    blackDot/.style={circle, draw, fill, inner sep=1.5pt},
    blueDot/.style={circle, draw, blue, fill, inner sep=1.5pt},
    cyanDot/.style={circle, draw, cyan, fill, inner sep=1.5pt},
    redEdge/.style={red, line width=0.5mm},
    blackEdge/.style={black, line width=0.5mm},
    transp/.style={inner sep=1.5pt},
    redDotted/.style={red, densely dotted, line width=0.5mm},
    blueEdge/.style={line width = 0.5mm, blue, ->, -stealth},
}

\makeatletter
\pgfdeclareshape{myNode}{
\inheritsavedanchors[from=rectangle] 
\inheritanchorborder[from=rectangle]
\inheritanchor[from=rectangle]{center}
\inheritanchor[from=rectangle]{north}
\inheritanchor[from=rectangle]{south}
\inheritanchor[from=rectangle]{west}
\inheritanchor[from=rectangle]{east}
\backgroundpath{
    \southwest \pgf@xa=\pgf@x \pgf@ya=\pgf@y
    \northeast \pgf@xb=\pgf@x \pgf@yb=\pgf@y
    \pgfsetcornersarced{\pgfpoint{5pt}{5pt}}
    \pgfpathmoveto{\pgfpoint{\pgf@xa}{\pgf@ya}}
    \pgfpathlineto{\pgfpoint{\pgf@xa}{\pgf@yb}}
    \pgfpathlineto{\pgfpoint{\pgf@xb}{\pgf@yb}}
    \pgfsetcornersarced{\pgfpoint{5pt}{5pt}}
    \pgfpathlineto{\pgfpoint{\pgf@xb}{\pgf@ya}}
    \pgfpathclose
}
}
\makeatother

\usepackage{xcolor}

\def\edgeSizeLB{1}
\def\depthFig{-1.5}
\def\depthOneX{2}
\def\depthOneY{\depthFig}

\def\depthTwoY{\depthOneY + \depthFig}
\def\depthThreeX{2}
\def\depthThreeY{\depthTwoY + \depthFig}
\def\edgeSizeLB{1}
\def\depthFig{-1.5}
\def\depthOneX{2}
\def\depthOneY{\depthFig}

\def\depthTwoY{\depthOneY + \depthFig}
\def\depthThreeX{2}
\def\depthThreeY{\depthTwoY + \depthFig}

\title[Beyond the Pseudoforest Strong Nine Dragon Tree Theorem]{Beyond the Pseudoforest Strong Nine Dragon Tree Theorem}

\author[Mies]{Sebastian Mies}
\address[Sebastian Mies]{Institute of Computer Science, Johannes Gutenberg University Mainz}
\email{smies@students.uni-mainz.de}
\author[Moore]{Benjamin Moore} \thanks{This work was completed while Benjamin Moore was a postdoc at Charles University, supported by  project 22-17398S (Flows and cycles in graphs on surfaces) of Czech Science Foundation.}
\address[Benjamin Moore]{Institute of Science and Technology, Klosterneuburg, Austria }   
\email{Benjamin.Moore@ist.ac.at}
\author[Smith-Roberge]{Evelyne Smith-Roberge}
\address[Evelyne Smith-Roberge]{School of Mathematics, Georgia Institute of Technology}
\email{esmithroberge3@gatech.edu}

\date{}

\newtheorem{thm}{Theorem}[section]

\newtheorem{lemma}[thm]{Lemma}

\newtheorem{conj}[thm]{Conjecture}
\newtheorem{cor}[thm]{Corollary}
\newtheorem{claim}{Claim}
\newtheorem{question}[thm]{Question}

\theoremstyle{definition}
\newtheorem{definition}[thm]{Definition}
\newtheorem{obs}[thm]{Observation}

\newtheorem*{ack}{Acknowledgements}

\newtheoremstyle{case}{}{}{\normalfont}{}{\itshape}{\normalfont:}{ }{}

\theoremstyle{case}

\numberwithin{equation}{section}

\newcommand{\floor}[1]{\left\lfloor#1\right\rfloor }
\newcommand{\ceil}[1]{\left\lceil#1\right\rceil }

\date{}

\begin{document}
\maketitle

\def\lowerMadBoundWithAlpha{\frac{\ell(k+1) + \alpha}{(\ell + 1)(k+1) + \alpha}}
\def\lowerMadBoundWithOne{\frac{(k+1)\ell + 1}{(k+1)(\ell + 1) + 1}}
\def\lowerMadBoundWithZero{\frac{(k+1)\ell}{(k+1)(\ell + 1)}}
\def\seb#1{{\color{green} #1}}
\def\sebsout#1{{\color{green} \sout{#1}}}
\def\sebcancel#1{{\color{green} \cancel{#1}}}
\def\markRed#1{{\color{red} #1}}
\def\kAryToPath#1{{\color{orange} #1}}
\def\kAryToPathSout#1{\kAryToPath{\sout{#1}}}
\def\kAryToPathCancel#1{\kAryToPath{\cancel{#1}}}
\def\compSizeBoundWithZ#1{d - #1(k-1) + 1}
\def\KC{K_{\mathcal C}}
\def\densOfKC{\frac{e(\KC)}{v(\KC)}}
\def\densNDT{\frac{d}{d+k+1}}
\def\equivModKPlusOne#1#2{#1 \equiv #2 \mod (k+1)}
\def\nequivModKPlusOne#1#2{#1 \not\equiv #2 \mod (k+1)}
\def\fracArb#1{\gamma(#1)}
\def\explSG{H_{f, R}}

\begin{abstract}
    The pseudoforest version of the Strong Nine Dragon Tree Conjecture states that if a graph $G$ has maximum average degree $mad(G) = 2 \max_{H \subseteq G} \frac{e(H)}{v(H)}$ at most $2(k + \densNDT)$, then it has a decomposition into $k+1$ pseudoforests where in one pseudoforest $F$ the components of $F$ have at most $d$ edges. This was proven in 2020 in \cite{sndtcPsfs}. We strengthen this theorem by showing that we can find such a decomposition where additionally $F$ is acyclic, the diameter of the components of $F$ is at most $2\ell + 2$, where $\ell = \floor{\frac{d-1}{k+1}}$, and at most $2\ell + 1$ if $\equivModKPlusOne{d}{1}$. Furthermore, for any component $K$ of $F$ and any $z \in \mathbb N$, we have $diam(K) \leq 2z$ if $e(K) \geq d - z(k-1) + 1$. We also show that both diameter bounds are best possible as an extension for both the Strong Nine Dragon Tree Conjecture for pseudoforests and its original conjecture for forests. In fact, they are still optimal even if we only enforce $F$ to have any constant maximum degree, instead of enforcing every component of $F$ to have at most $d$ edges. 
\end{abstract}
\section{Introduction}
In this paper we study decompositions. Recall that a \textit{decomposition} of a graph $G$ is a partition of the edge set of the graph. We will use the notation that $e(G)$ is the number of edges of $G$, and $v(G)$ is the number of vertices of $G$. Decompositions of graphs are a heavily studied area, as if one can decompose a graph into only a few simple pieces, then one can in some sense deduce that the entire graph structure is simple.  One of the most simple classes of graphs is the class of  forests; thus a natural question is \textit{when does a graph decompose into $k$ forests?}  This question has been answered completely by a beautiful theorem of Nash-Williams:

\begin{thm}[\cite{nash}]\label{nashthm}
A graph $G$ decomposes into $k$ forests if and only if 
\[ \fracArb{G} := \max_{\substack{H \subseteq G,\\ v(H) \geq 2}} \frac{e(H)}{v(H)-1} \leq k.\]
\end{thm}

We will refer to $\fracArb{G}$ as the \textit{fractional arboricity of $G$}. The notation $H \subseteq G$ means that $H$ is a subgraph of $G$. Noting that a forest $F \subseteq G$ has at most $v(G) -1$ edges, and that if a graph $G$ decomposes into $k$ forests, then so do all of its subgraphs, we see that Nash-Williams' Theorem says that the obvious necessary conditions are in fact sufficient. 

A natural approach to strengthening Nash-Williams' Theorem would be to try and control the types of forests that appear in the decomposition. In general, not much can be done, but by observing that the fractional arboricity need not be integral, one might intuitively believe that if a graph $G$ has fractional arboricity closer to $k$ than to $k+1$, then $G$ is very close to decomposing into $k$ forests instead of $k+1$ forests, and thus one could impose more structure on at least one of the forests in the decomposition. For example, applying Nash-Williams' Theorem to cycles gives us that cycles decompose into two forests as they have fractional arboricity exactly $\frac{n}{n-1}$ if $n$ is the number of vertices in the cycle. However, of course a cycle decomposes into a forest and a matching, which is a significantly stronger statement than a decomposition into two forests.  More substantially, let us consider what happens for planar graphs. By Euler's formula, a simple planar graph $G$ with at least three vertices has at most $3v(G)-6$ edges, and if $G$ is moreover triangle-free, then $G$ has at most $2v(G) -4$ edges. Thus Nash-Williams' Theorem says that simple planar graphs decompose into at most $3$ forests, and triangle-free planar graphs decompose into at most $2$ forests. Assuming the planar graph is not a forest, Nash-Williams' Theorem cannot say anything further even if $G$ has larger girth. But it is easy to see that the average degree of a simple planar graph of girth $g$ is bounded by $\frac{2g}{g-2}$ (see \cite{dischargingtutorial}), and thus the fractional arboricity of planar graphs tends towards $1$ as the girth increases. Multiple papers have been written showing that better decompositions can be found for planar graphs of larger girth: two of the authors of this paper showed that planar graphs with girth at least $5$ decompose into two forests, one of which has every component containing at most five edges \cite{miesMoore}; Kim et al. \cite{Kostochkaetal} showed that planar graphs of girth at least $6$ decompose into a forest and a forest where each component has at most $2$ edges; and Montassier et al. \cite{ndtConjs} showed that planar graphs of girth at least $8$ decompose into a forest and a matching. 
The Strong Nine Dragon Tree Conjecture, proposed in \cite{ndtConjs}, formalizes this intuition:

\begin{conj}[Strong Nine Dragon Tree Conjecture]
For any integers $k$ and $d$, every graph $G$ with fractional arboricity at most $k+ \frac{d}{k+d+1}$ decomposes into $k+1$ forests such that in one of the forests, every connected component contains at most $d$ edges.
\end{conj}

This conjecture remains wide open despite a large effort. In \cite{ndtConjs}, Montassier et al.  proved the $k=1$ and $d=1$ case. In \cite{Kostochkaetal}, Kim et al. proved the $k=1$ and $d=2$ case. In a breakthrough result, Yang proved the Strong Nine Dragon Tree conjecture when $d=1$ and for every $k$ \cite{Yangmatching}. 
In \cite{approxArb} Blumenstock and Fischer proved the $k=d$ case for some special graph classes. Recently, two of the authors proved the conjecture when $d \leq k+1$ and gave an approximate bound when $d \leq 2(k+1)$ (see \cite{miesMoore}). The strongest evidence for the Strong Nine Dragon Tree Conjecture comes from the Nine Dragon Tree Theorem, proven by Jiang and Yang:

\begin{thm}[Nine Dragon Tree Theorem, \cite{ndtt}]
For any integers $k$ and $d$, every graph $G$ with fractional arboricity at most $k+ \frac{d}{k+d+1}$ decomposes into $k+1$ forests, such that one of the forests has maximum degree at most $d$.
\end{thm}

In fact, the Nine Dragon Tree Theorem is tight in the following sense:

\begin{thm}[\cite{ndtConjs}]
\label{montassierexample}
 For any integers $k$ and $d$, there exist infinitely many graphs $G$ such that there exists a set  $D \subseteq E(G)$ where $|D| \leq d$, $\fracArb{G-D} = k + \frac{d}{k+d+1}$, but $G$ does not decompose into $k+1$ forests where one of them has maximum degree $d$.
\end{thm}

The idea of the construction given in Theorem \ref{montassierexample} is that there exists a graph $G-D$ where in every decomposition of $G$ into $k+1$ forests where one of the forests has maximum degree at most $d$, every component of this forest is an isolated vertex, or a star with exactly $d$ edges. By then cleverly adding few edges, one can increase the fractional arboricity only very slightly, while preventing a decomposition that satisfies the conclusion of the Nine Dragon Tree Theorem. In light of this, Xuding Zhu asked if the Strong Nine Dragon Tree Conjecture could be strengthened to the following:

\begin{question}[Zhu, personal communication] \label{conj:zhuStars}
 For any integers $k$ and $d$, does every graph $G$ with fractional arboricity at most $k+ \frac{d}{k+d+1}$ decompose into $k+1$ forests such that the components of one of the forests are isomorphic to stars each containing at most $d$ edges?
\end{question}

We will unfortunately show the answer to this question is no and we conjecture that it is still no even if we remove the constraint on the number of edges of the components. Had the answer been yes, this would have had applications that a proof of the Strong Nine Dragon Tree Conjecture would not. For example, Merker and Postle \cite{boundeddiameter} showed that if a graph decomposes into a forest and a forest where all components are isomorphic to stars (called star forests), then the  graph decomposes into two forests, where every component has diameter\footnote{Recall the \textit{distance} between vertices $u$ and $v$ in a graph $G$ is the minimum number of edges in a path with endpoints $u$ and $v$. The \textit{diameter} of $G$ is the maximum distance between any two vertices in $G$.} at most 18.   So had the answer to Question \ref{conj:zhuStars} been yes, using Euler's formula, it would have implied that planar graphs of girth at least five decompose into a forest and a star forest, and thus into two forests where each component has diameter at most $18$ (and in fact $14$ from an improvement in \cite{mousaviHaji}).

Even though the Strong Nine Dragon Tree Conjecture remains open, the Nine Dragon Tree Theorem as well as the partial progress theretoward  have already had applications. For example, the Nine Dragon Tree Theorem as well as a result of Zhu \cite{ZHUgamecol} imply best possible bounds on the game chromatic number of planar graphs with girth at least five, and two of the authors used the partial progress towards the Strong Nine Dragon Tree Conjecture to show that all $5$-edge-connected planar graphs have a $\frac{5}{6}$-thin spanning tree \cite{miesMoore}. 

While the original motivation for the Nine Dragon Tree Conjectures came from its applications, there is no reason to restrict our attention to just forests. One can pose a Strong Nine Dragon Tree-style conjecture for any family of graphs that has a Nash-Williams type theorem. In fact, there are pseudoforest versions \cite{sndtcPsfs, ndttPsfs}, digraph versions \cite{digraphndt}, and even matroidal versions \cite{matroidndt}. For this paper, we will be interested in \textit{pseudoforests}. Recall that a pseudoforest is a graph where each connected component contains at most one cycle. Hakimi showed:

\begin{thm}[\cite{hakimi}, Hakimi's Theorem]
A graph $G$ decomposes into $k$ pseudoforests if and only if 
\[\text{mad}(G) := \max_{\substack{H \subseteq G, \\ V(H) \neq \varnothing}} \frac{2e(H)}{v(H)} \leq 2k\]
\end{thm}

Again, just like Nash-Williams' Theorem, Hakimi's Theorem says that the obvious necessary condition is sufficient. We will refer to $mad(G)$  as the \textit{maximum average degree} of a graph $G$. Similarly to fractional arboricity, the maximum average degree of a graph need not be integral. Fan et al. \cite{ndttPsfs} proved a pseudoforest version of the Nine Dragon Tree Conjecture:

\begin{thm}[Pseudoforest Nine Dragon Tree Theorem]
Let $k$ and $d$ be integers. Every graph with maximum average degree at most $2(k + \densNDT)$ decomposes into $k+1$ pseudoforests where one of the pseudoforests has maximum degree $d$. Further, for every $\epsilon >0$ and every pair $k$ and $d$, there exists a graph $G$ such that $G$ has maximum average degree at most $2(k+ \frac{d}{k+d+1} + \epsilon)$ but does not decompose into $k+1$ pseudoforests where one of the pseudoforests has maximum degree $d$.
\end{thm}

Similarly, we can ask for a pseudoforest analogue of the Strong Nine Dragon Tree Conjecture. This was proven by Grout and one of the authors \cite{sndtcPsfs}:

\begin{thm}[Pseudoforest Strong Nine Dragon Tree Theorem]
\label{thm:psndt}
Let $k$ and $d$ be integers. Every graph with maximum average degree at most $2(k + \frac{d}{k+d+1})$ decomposes into $k+1$ pseudoforests where in one of the pseudoforests, every connected component contains at most $d$ edges.
\end{thm}

In this paper, we strengthen this theorem, showing that not only can we make the components of one of the pseudoforests have at most $d$ edges, but that we can make this pseudoforest a forest, and further bound the diameter of every component.

\begin{restatable}{thm}{restupperbound}
\label{thm:upperBound}%
Let $k, d \in \mathbb N$, where $k \geq 1$. Let $\ell = \floor{\frac{d-1}{k+1}}$. Then every graph $G$ with maximum average degree at most $2(k + \frac{d}{d+k+1})$ decomposes into  $k + 1$ pseudoforests where one of the pseudoforests $F$ satisfies the following:
\begin{itemize}
\item $F$ is acyclic,
\item every component $K$ of $F$ has $e(K) \leq d$,
\item $diam(K) \leq 2\ell + 2$, and if $\equivModKPlusOne{d}{1}$, then $diam(K) \leq 2\ell+1$,
\item for every component $K$ of $F$ satisfying $e(K) \geq d -z(k-1) +1$, we have $diam(K) \leq 2z$ for any $z \in \mathbb{N}$.
\end{itemize}
\end{restatable}

We remark that Theorem \ref{thm:upperBound} holds even if the graphs are allowed to contain loops and parallel edges.
Thus we show that in some sense, the spirit of Question \ref{conj:zhuStars} is maintained for pseudoforests, even though we cannot always make the last pseudoforest a star forest where all of the components have at most $d$ edges. We point out an appealing corollary where this does in fact hold:

\begin{cor}
\label{cor:starforests}
Let $k$ and $d$ be integers such that $d \leq k+1$. Every graph with maximum average degree at most $2(k + \frac{d}{d+k+1})$ decomposes into $k + 1$ pseudoforests where one of the pseudoforests $F$ is a forest, every component $K$ of $F$ has $e(K) \leq d$ and is isomorphic to a star. 
\end{cor}

A natural class where one can apply Corollary \ref{cor:starforests} is the class of graphs with fixed odd maximum degree. Let $G$ be a graph with odd maximum degree $\Delta$. It follows that $G$ has maximum average degree at most $\Delta$. Setting $k = \frac{\Delta-1}{2}$, and $d = \frac{\Delta+1}{2}$, gives that $2(k + \frac{d}{k+d+1}) = \Delta$. Further, $d \leq k+1$. Therefore such graphs decompose into $\frac{\Delta-1}{2} +1$ pseudoforests, where one of the pseudoforests is a star forest where every component has at most $\frac{\Delta+1}{2}$ edges.

In addition, we construct examples showing our theorem is best possible in the following sense.

\begin{restatable}{thm}{restlowerBound}
\label{thm:lowerBound}%
 Let $k, \ell, D \in \mathbb N$, $\epsilon > 0$, $k \geq 1$ and $\alpha \in \{0, 1\}$.
    There are simple graphs $G$ with $\fracArb{G} < k + \lowerMadBoundWithAlpha + \epsilon$ that do not have a decomposition into $k+1$ pseudoforests where one of the pseudoforests has maximum degree at most $D$ and the diameter of every component is less than $2\ell + 1 + \alpha$.
\end{restatable}

Note that for $\ell = \floor{\frac{d-1}{k+1}}$ we have $\lowerMadBoundWithZero < \densNDT$ and if $\equivModKPlusOne{d}{1}$, then $\lowerMadBoundWithOne = \densNDT$. Furthermore, our second statement that bounds the diameter of large components even more if $k \geq 2$ is also best possible in a similar way.

\begin{restatable}{thm}{restzlowerBound}
\label{thm:z:lowerBound}
 Let $k, D, d, z \in \mathbb N$, $\epsilon > 0$, $k \geq 2$, $z \leq \floor{\frac{d-1}{k+1}}$. There are simple graphs $G$ with $\fracArb{G} < k + \densNDT + \epsilon$, where in any decomposition into $k+1$ pseudoforests such that one pseudoforest $F$ has maximum degree at most $D$, $F$ has a component $K$ with $e(K) \geq d - z(k-1) + 1$ and $diam(K) \geq 2(z+1)$.
\end{restatable}

Note that since $mad(G)/2 < \fracArb{G}$ we can replace $\fracArb{G}$ with $mad(G)/2$ in both Theorem \ref{thm:lowerBound} and \ref{thm:z:lowerBound} to prove the optimality of Theorem \ref{thm:upperBound}. On the other hand we can replace the words \textit{pseudoforest} with \textit{forest} in order to get a lower bound for diameter refinements of the Strong Nine Dragon Tree Conjecture.

We now sketch how we prove these theorems. For Theorem \ref{thm:upperBound}, we follow the approach first developed for pseudoforests in \cite{ndttPsfs} and extended to some forests in \cite{Yangmatching}; which later was pushed further to prove the Nine Dragon Tree Theorem in \cite{ndtt}, then the Pseudoforest Strong Nine Dragon Tree Theorem in \cite{sndtcPsfs}; and which was finally used to give the the current best known result on the Strong Nine Dragon Tree Conjecture \cite{miesMoore}. Although focusing on establishing the diameter constraints, we will reprove Theorem \ref{thm:psndt} without much additional effort. 
 
The idea of the proof of Theorem \ref{thm:upperBound} is as follows: we start with some decomposition into $k+1$ pseudoforests and we pick any of the pseudoforests, say $F$, and try to transform $F$ into a forest with bounded component size and diameter. We may assume that the remaining $k$ pseudoforests are maximal (i.e., have $v(G)$ edges), as otherwise we can rearrange the decomposition to reduce the number of edges in $F$. Furthermore, we may assume that $F$ either contains a cycle, a component with too many edges, or a component with too large diameter; such components we will call \textit{bad}. We will search around this bad component to try to find nearby components that are small and acyclic which we can use to ``augment" our decomposition. To formalize this, we will use the well-known fact that a graph decomposes into $k$ pseudoforests if and only if it has an orientation with outdegree at most $k$. Using the orientation of the other $k$ pseudoforests, we will be able to prove that if a bad component is close to many  acyclic components with few edges (from here on out, called \textit{small} components), we can augment our decomposition such that we can remove the bad component, while maintaining a pseudoforest decomposition. Unfortunately, there is no guarantee that a bad component will be near many small components and so we will explore further from the components near the bad component to try and find an augmentation which will make the components near the bad component small, so that we can then augment and get rid of the bad component. Continuing the exploration until we can no longer do so, we get what we call an \textit{exploration subgraph}. In this exploration subgraph, we put an order on the components (which we call a \textit{legal order}) that encodes in some sense both the distance of components from the bad component, as well as their sizes. We will show that in this subgraph, either we can augment the decomposition to make components closer to the bad component smaller or acyclic (which we call \textit{reducing} the legal order), which would eventually result in many components near the bad component being small (and thus allowing us to augment our decomposition and get rid of the bad component), or we cannot do any of these augmentations. In this latter case,  we will show this implies that the exploration subgraph actually contains a large number of edges, contradicting the maximum average degree bound. The augmentations performed are similar to those in the proof of Theorem \ref{thm:psndt}, but we require a much more refined analysis to obtain Theorem \ref{thm:upperBound}. 

For Theorems \ref{thm:lowerBound} and \ref{thm:z:lowerBound} we do the following: we start with an exploration subgraph $C$ in which $F$ contains a component with large diameter and which cannot be improved any further by our algorithm for the upper bounds. We then add a set $S$ of just a few vertices as well as many copies of $C$ adding edges between $S$ and $C$ carefully. Since each of the vertices of $S$ can only have a few adjacent edges in $F$, we can force the edges between $S$ and at least one copy of $C$ to belong to the other $k$ pseudoforests. This will then force many edges within $C$ to belong to $F$, which results in $F$ having long paths.

We end the introduction with some notation and conventions.  All mentioned paths are meant to be simple. If $G$ is a graph and $X \subseteq V(G)$, then $G[X]$ denotes the \textit{induced subgraph} of $G$ by $X$. We also write $E(X)$ for $E(G[X])$ and $e(X)$ for $e(G[X])$ if it is clear from context which graph $G$ is underlying. Further, $deg_X(v)$ denotes the degree and $deg^+_X(v)$ the outdegree of a vertex $v \in X$ in $G[X]$. If $X = V(G)$ and it is clear which underlying graph $G$ we have, we simply write $deg(v)$ or $deg^+(v)$. If $A$ is a set, then let $A + a := A \cup \{a\}$ and $A - a := A \setminus \{a\}$. If $G$ is a graph and $S \subseteq V(G)$, then $G - S$ is the graph obtained by removing the vertices of $S$ and all incident edges from $V(G)$ and $E(G)$, respectively. Further, let $G - v := G - \{v\}$ for vertices $v$ of $G$. Similarly, if $e$ is an edge of a graph $G$, then $G - e := (V(G), E(G) - e)$. Let $E(V_1, V_2)$ denote the set of all edges in $G$ with exactly one endpoint in each of the sets $V_1, V_2 \subseteq V(G)$, and let $e(V_1, V_2) := |E(V_1, V_2)|$. 

The structure of the paper is as follows: in Section 2 we will prove Theorem \ref{thm:upperBound} and in Sections 3 and 4 we prove Theorems \ref{thm:lowerBound} and \ref{thm:z:lowerBound}, respectively.

\section{Proof of the upper bounds}

Our goal in this section is to prove Theorem \ref{thm:upperBound}, restated below for convenience.

\restupperbound*

Throughout this section, graphs are allowed to contain parallel edges and loops. Before we begin, we mention a well-known correspondence between pseudoforests and orientations that we will exploit. Recall that, given a graph $G$, an \textit{orientation} of $G$ is obtained by taking each edge $xy \in E(G)$, and replacing $xy$ with exactly one of the arcs $(x,y)$ and $(y,x)$.

\begin{obs}
A graph $G$ is a pseudoforest if and only if $G$ admits an orientation where every vertex has outdegree at most one. 
\end{obs}
From this observation, we get an important corollary.

\begin{cor}\label{cor:pseudodecompiff}
A graph admits a decomposition into $k$ pseudoforests if and only if it admits an orientation such that every vertex has outdegree at most $k$.
\end{cor}

\subsection{Picking the counterexample}

We prove Theorem \ref{thm:upperBound} by contradiction. In this subsection, we describe how we will pick our counterexample and explain the setup for our proof. To that end: let $k$ and $d$ be fixed positive integers, and suppose that $G$ is a vertex-minimal counterexample to Theorem \ref{thm:upperBound} for the{se} fixed values of $k$ and $d$. 

Our first step will be to obtain desirable orientations of $G$. For this, we use a lemma proved in \cite{ndttPsfs} (Lemma $2.1$). Technically, we need a stronger lemma, but the same proof as Lemma $2.1$ in \cite{ndttPsfs} suffices to prove the strengthening given below.

\begin{lemma}[\cite{ndttPsfs}]
    \label{orientationlemma}
    If $G$ is a vertex-minimal counterexample to Theorem \ref{thm:upperBound}, then there exists an orientation of $G$ such that for all $v \in V(G)$, we have $k \leq deg^+(v) \leq k+1$. 
\end{lemma}

Let $\mathcal{F}$ be the set of edge-coloured (mixed) subgraphs of $G$ where every vertex $v \in V(G)$ has exactly $k$ outgoing blue (directed) edges and all the remaining edges are red, undirected, and induce a pseudoforest. We will call $f \in \mathcal F$ a \textit{red-blue colouring}. Note that $\mathcal F$ is non-empty by Lemma \ref{orientationlemma} and Corollary \ref{cor:pseudodecompiff}. 
Furthermore, we can extract a decomposition into $k+1$ pseudoforests from $f$, where $k$ pseudoforests are coloured blue and one is coloured red. If $k \geq 2$, we have multiple options to form the blue pseudoforests. The red pseudoforest will be the one in which we wish to eliminate cycles and bound the size and diameter of each connected component. As a note for the reader, as we modify a red-blue colouring in the following pages, it will be easier to focus on maintaining a blue $k$-orientation rather than considering the $k$ individual blue pseudoforests.
Furthermore, note that we call the components of a red pseudoforest \textit{red components}.\\
We now categorize the red components that contradict Theorem \ref{thm:upperBound}.

\begin{definition} \label{def:bad}
    A red component $K$ of a pseudoforest $F$ is \textit{bad} if $K$ satisfies any of the following: 
    \begin{enumerate}
        \item contains a cycle, or
        \item has more than $d$ edges, or
        \item has $diam(K) > 2\ell + 2$, or
        \item has $diam(K) = 2\ell + 2$ when $d \equiv 1 \mod (k+1)$, or
        \item has $diam(K) > 2z$ and $e(K) \geq \compSizeBoundWithZ{z}$ for a $z \in \mathbb N$ with $1 \leq z \leq \ell$.
    \end{enumerate}
    We say $K$ is \textit{(i)-bad} if $(i)$ is true for $K$ and $(j)$ is not true for $K$ for all $j < i$.
\end{definition}

Note that for (5)-bad components we chose the range $1 \leq z \leq \ell$ since if $z = 0$ then $K$ is (2)-bad, and if $z > \ell$ then $K$ is (3)-bad.

As $G$ is a counterexample, in every red-blue colouring there is at least one bad component (of the red pseudoforest). Our goal will be to augment the decomposition in order to reduce the number of bad components, or possibly make a bad component ``less bad". We will not be able to do that in a single step necessarily,  so instead we will focus on one bad component and a specific subgraph that stems from this component where we will perform possibly many augmentations to eventually make the bad component less bad. To that end, we define the following subgraph.
\begin{definition}
    Suppose that $f$ is a red-blue colouring of $G$. Let $R$ be a bad component of this colouring. We define the \textit{exploration subgraph $\explSG$ of $f$ rooted at $R$} in the following manner: let $S \subseteq V(G)$ where $v \in S$ if and only if there exists a path $P = v_{1},\ldots,v_{m}$ such that $v_{m}= v$, $v_{1} \in V(R)$, and for each $i \in \{1, \dots, m-1\}$ either $(v_{i},v_{i+1})$ is a blue arc or $v_{i}v_{i+1}$ is a red edge. Let $\explSG$ be the graph induced by $S$, where $E(\explSG)$ inherits the colours of $f$ as well as the orientations of the blue arcs (whereas the red edges of $\explSG$ remain undirected).
\end{definition}

It might not yet be clear why we made this particular definition for $\explSG$; however, the next observation shows that for any exploration subgraph $\explSG$, the red edge density must be low. Before stating this observation, we fix some notation. Given a subgraph $G'$ of $G$, we will let $E_{b}(G')$ and $E_{r}(G')$ denote the sets of edges of $G'$ coloured blue and red, respectively. We let $e_{b}(G') = |E_{b}(G')|$ and $e_{r}(G') = |E_{r}(G')|$.

\begin{obs}
    \label{finishingobservation}
    For any red-blue colouring $f$ of $G$ and any choice of root component $R$, the exploration subgraph $\explSG$ satisfies
    \[\frac{e_{r}(\explSG)}{v(\explSG)} \leq \densNDT.\]
\end{obs}

\begin{proof}
    Suppose towards a contradiction that $e_{r}(\explSG) / v(\explSG) > d / (d+k+1).$
    As $\explSG$ is an induced subgraph defined in part by directed paths and every vertex $v \in V(G)$ has $k$ outgoing blue edges by definition of $f$, each vertex in $\explSG$ has $k$ outgoing blue edges.
    Thus, $e_b(\explSG) / v(\explSG) = k.$ Therefore
    
    \[\frac{\text{mad}(G)}{2} \geq \frac{e(\explSG)}{v(\explSG)} = \frac{e_{r}(\explSG)}{v(\explSG)} + \frac{e_{b}(\explSG)}{v(\explSG)} > k + \densNDT.\]
    
    But this contradicts that $G$ has $\text{mad}(G) \leq 2(k+\densNDT)$. 
    
\end{proof}

Throughout the proof of Theorem \ref{thm:upperBound}, we will be attempting to show that we can augment a given decomposition in such a way that either we obtain a decomposition satisfying Theorem \ref{thm:upperBound} or we can find an exploration subgraph $\explSG$ which contradicts Observation \ref{finishingobservation}.

As Observation \ref{finishingobservation} allows us to focus only on red edges, it is natural to focus on red components which have small average degree. With this in mind, we define the notion of a \textit{small red component}.

\begin{definition}
    Let $K$ be a connected acyclic red component, and let $\ell$ be as in Theorem \ref{thm:upperBound}. We say $K$ is \textit{small} if  $e(K) \leq \ell$.
\end{definition}

The following lemma motivates the definition of $\ell := \floor{\frac{d-1}{k+1}}$.

\begin{lemma} \label{lemma:ellAndDensity}
    A red component $K$ is small if and only if
    $\frac{e(K)}{v(K)} < \densNDT$.
    In particular, $\ell$ is the largest integer $n$ such that $\frac{n}{n+1} < \densNDT$.
\end{lemma}
\begin{proof}
    If $K$ contains a cycle, then $K$ is not small and we have $\frac{e(K)}{v(K)} = 1 > \densNDT$.
    If $K$ is acyclic, $\frac{e(K)}{v(K)} \!=\! \frac{e(K)}{e(K) + 1} < \densNDT$ is equivalent to $e(K) < \frac{d}{k+1}$, which is equivalent to $e(K) \leq \floor{\frac{d-1}{k+1}} = \ell$.
\end{proof}

A small component cannot be bad; we prove this in Corollary \ref{cor:smallIsNotBad}.   We will want to augment our decomposition, and we will want a measure of progress that our decomposition is improving. Of course, reducing the number of bad components or their size would be a clear improvement. However, this might not always be possible, so we will introduce a notion of a ``legal order" of the red components. 
This order keeps track of how ``close'' a red component is to the root component. Our goal is to continually perform augmentations that make components ``closer" to the root component have fewer edges, without creating any new bad components. If our augmentations result in small red components next to the root component, then we might be able to perform an augmentation which makes the root component ``less bad'' or even not bad at all.

We formalize this in the following manner.  

\begin{definition}
    We call an order $(R_{1},\ldots,R_{t})$ of the red components of $\explSG$ \textit{legal} if all red components of $\explSG$ are in the order, $R_{1}$ is the root component $R$, and for each $j \in \{2,\ldots,t\}$ there exists an integer $i$ with $1 \leq i < j$ such that there is a blue arc $(u,v)$ with $u \in V(R_{i})$ and $v \in V(R_{j})$.
\end{definition} 

Let $(R_{1},\ldots,R_{t})$ be a legal order. We will say that $R_{i}$ is a \textit{parent} of $R_{j}$ if $i < j$ and there is a blue arc $(v_{i},v_{j})$ in $\explSG$ where $v_{i} \in V(R_{i})$  and $v_{j} \in V(R_{j})$. Note that a red component may have many parents.

If $R_{i}$ is a parent of $R_{j}$, then we say that $R_{j}$ is a \textit{child of $R_{i}$ generated by $(v_i, v_j)$}. We say a red component $R_{i}$ is an \textit{ancestor} of $R_{j}$ if there exists a sequence of red components $R_{i_{1}},\ldots,R_{i_{m}}$ such that $R_{i_{1}} = R_{i}$, $R_{j_{m}} = R_{j}$, and $R_{i_{q}}$ is a parent of $R_{i_{q+1}}$ for all $q \in \{1,\ldots,m-1\}$. An important concept is that of vertices \textit{witnessing} a legal order.

\begin{definition}
    Given a legal order $(R_{1},\ldots,R_{t})$ and integer $j \in \{2, \dots, t\}$, we say a vertex $v \in V(R_j)$ \textit{witnesses the legal order for $R_{j}$} if there is a blue arc $(u,v)$ and integer $1 \leq i< j$ such that $u \in V(R_{i})$. 
\end{definition}

Observe that there may be many vertices that witness the legal order for a given red component. More importantly, for every red component except the root, there exists a vertex that witnesses the legal order. We also want to compare two different legal orders. 

\begin{definition}
    Let $a = (a_{1},\ldots,a_{t})$ and $a' = (a'_{1},\ldots,a'_{t'})$ be two tuples of natural numbers. Let $a_i = 0$ for $t < i \leq t'$ and $a'_i = 0$ for $t' < i \leq t$. We will say $a$ is \textit{lexicographically smaller} than $a'$ if there is an $m \in \mathbb N$ such that $a_i = a'_i$ for all $i < m$ and $a_m < a'_m$. In this case we write $a <_{lex} a'$. Note that this defines a total order $\leq_{lex}$.
\end{definition}

We are now able to define a rating for red-blue colourings which indicates how close such a colouring is to satisfying Theorem \ref{thm:upperBound}. Namely: we say that an $(i)$-bad component is always worse than a $(j)$-bad component if $i < j$; that an $(i)$-bad component $K$ is worse than an $(i)$-bad component $K'$ if $e(K) > e(K')$; and if two colourings (with minimal legal orders $(R_1, \dots, R_t$) and $(R'_1, \dots, R'_{t'})$, respectively) have the same types of bad components, then we compare $(e(R_1), \dots, e(R_t))$ and $(e(R'_1), \dots, e(R'_{t'}))$ under $\leq_{lex}$; the idea is that each such tuple indicates how close the corresponding colouring is to an augmentation that reduces the ``badness'' of bad components. We will formalize this now.
Given a red-blue colouring $f$, let $\Delta^{(i)}(f) := \big(\Delta^{(i)}_{e(G)}(f), \Delta^{(i)}_{e(G)-1}(f), \ldots, \Delta^{(i)}_0(f)\big)$ where $\Delta^{(i)}_j(f)$ is the number of red components of $f$ that are $(i)$-bad and have exactly $j$ edges. Let $\sigma(f, R)$ denote the lexicographically smallest  tuple $(e(R_1), \ldots, e(R_t))$ over all legal orders $(R_1, \ldots, R_t)$ of $\explSG$. We call a legal order $(R'_1, \ldots, R'_{t'})$ 
 with $\sigma(f, R) = (e(R'_1), \ldots, e(R'_{t'}))$ a \textit{smallest legal order}.\\
For our vertex-minimal counterexample, we choose a colouring $f \in \mathcal F$ and a bad component $R$ of the red pseudoforest $F$ of $f$ such that $\Delta(f, R) = \big(\Delta^{(1)}(f), \ldots, \Delta^{(5)}(f), \sigma(f, R)\big)$ is lexicographically smallest.\\
For the occasions where we only want to compare the bad components of two red pseudoforests, we define $\Delta'(f) = \big(\Delta^{(1)}(f), \ldots, \Delta^{(5)}(f)\big)$.

From here on out, we will fix $f, R, F$ picked in the manner described above. If not stated otherwise, parental relationships between red components will be considered in the context of an arbitrary but fixed smallest legal order $(R_1, \ldots, R_t)$.

In order to augment our decomposition we will only use the following simple procedure taken from \cite{sndtcPsfs}.

\begin{definition}
    Let $K$ be a red component of $\explSG$ and let $C$ be an acyclic child of $K$ generated by $(x, y)$. Suppose that $e = xv$ is a red edge in $K$ incident to $x$. To \textit{exchange $e$ and $(x,y)$} is to perform the following procedure: first, change the colour of $(x,y)$ to red and remove its orientation. Second, change the colour of $e$ to blue and orient it to $(x,v)$.
\end{definition}

\begin{obs}[\cite{sndtcPsfs}]
    \label{flipobservation}
    Suppose we exchange the edge $e = xv$ and $(x,y)$. Then the resulting red-blue colouring is in $\mathcal F$.
\end{obs}

 We will implicitly use Observation \ref{flipobservation} throughout the paper.

Eventually, we will need to show that if we cannot augment our decomposition then the average degree of $\explSG$ is too high, contradicting Observation \ref{finishingobservation}. 
It is easy to show that if $\explSG$ does not contain any small red components, then the average degree of $\explSG$ is too high, since $\frac{e(K)}{v(K)} \geq \densNDT$ if and only if $K$ is not small. However, getting rid of all small components is not realistic. To get around this, we partition the components of $F$ so that each part contains a non-small component along with some of its small children.
If the average degree of every part is at least $\densNDT$, then we still obtain the same contradiction. This partitioning requires that small components always have a parent that is not small. We will prove this later in Corollary \ref{smallchildrencorollary}.

\begin{definition}
    Denote the set of red components of $\explSG$ that are not small by $\mathcal K$. In an arbitrary fashion we assign each small component to exactly one of its parents in  $\mathcal K$.\\
    Let $K \in \mathcal K$ and let $C_1, \dots, C_q$ be the small children of $K$ that were assigned to $K$.\\
    Then
    $\mathcal{C}(K) := \{C_1, \dots, C_q\}$ and 
    \[\KC := \big(V(K) \cup \bigcup_{C \in \mathcal{C}(K)} V(C), \;E(K) \cup \bigcup_{C \in \mathcal{C}(K)} E(C)\big).\]
\end{definition}

\begin{lemma} \label{lemma:howToGetMADContradiction}
    Assume that in $\explSG$
    \begin{enumerate}
        \item small components do not have small children,
        \item we have $\densOfKC > \densNDT$ for every bad component $K$, and
        \item we have $\densOfKC \geq \densNDT$ for every non-bad red component $K$.
    \end{enumerate}

    Then we obtain a contradiction to Observation \ref{finishingobservation}, which proves Theorem \ref{thm:upperBound}.
\end{lemma}
\begin{proof}
As small components do not have small children by assumption, it follows that any red component of $\explSG$ is a subgraph of $\KC$ for a non-small component $K$.
We have that $\densOfKC \geq \densNDT$ holds for every $\KC$ of this partition of $(V(\explSG , E_r(\explSG))$. Since there is at least one (non-small) bad component for which this inequality is strict, we obtain a contradiction to Observation \ref{finishingobservation}.
\end{proof}

The setup for our proof is now done. It remains to show the three conditions of Lemma \ref{lemma:howToGetMADContradiction} hold for our optimally chosen colouring $f$. 

In the following subsections we will show the second condition for all types of bad components. In the final Subsection \ref{sec:bad5nonBad} we will also show the first and third condition.

\subsection{Density of \boldmath\texorpdfstring{$\KC$}{K\_C} if \boldmath\texorpdfstring{$K$}{K} is (1)-bad}

In this subsection, we show that if we have a component which is (1)-bad, then it has no small children, and thus contributes to showing that the average degree of our exploration subgraph is too high.

\begin{obs} \label{redcyclesaturationlemma}
    If $C$ is an acyclic child of $K$ generated by $(x,y)$, then $x$ does not lie in a cycle of $F$. 
\end{obs}
\begin{proof}
    Suppose towards a contradiction that $x$ lies in a cycle of $F$. Let $e$ be an edge incident to $x$ which lies in the cycle coloured red. Now exchange $e$ and $(x,y)$. As $(x,y)$ was an arc between two distinct red components and $e$ was in the cycle coloured red, after performing the exchange, we reduce the number of cycles in $F$ by one and do not affect other cyclic components. Thus the exchange results in a red-blue colouring with a smaller $\Delta^{(1)}$ which is a contradiction.
\end{proof}

\begin{lemma}   \label{lemma:cyclicDoesntHaveAcyclicChildren}
    If $K$ is a cyclic red component of $\explSG$, then there is no blue arc $(x, y)$ from $K$ to an acyclic red component $C$.
\end{lemma}
\begin{proof}
    Suppose towards a contradiction that there is such an edge. By Observation \ref{redcyclesaturationlemma} we know that $x$ does not lie on a cycle of $F$. There is a unique red path from $x$ to the cycle of $K$. On this path let $w$ be the neighbour of $x$. Exchange $xw$ and $(x, y)$.
    In the resulting red pseudoforest, $x$ and $y$ are in an acyclic component and $w$ is in a cyclic component containing fewer edges than $K$.
    Again, this results in a red-blue colouring with smaller $\Delta^{(1)}$, a contradiction.
\end{proof}

We are now equipped to prove the main result of this subsection.
\begin{lemma}
\label{fracbound1bad}
If a red component $K$ of $\explSG$ is (1)-bad, then $\densOfKC > \densNDT$.
\end{lemma}
\begin{proof}
    If $K$ is cyclic, then $K$ does not have any small children by Lemma \ref{lemma:cyclicDoesntHaveAcyclicChildren}. Thus, 
    \[
        \densOfKC = \frac{e(K)}{v(K)} = 1 > \densNDT.
    \]
\end{proof}

\subsection{Density of \boldmath\texorpdfstring{$\KC$}{K\_C} if \boldmath\texorpdfstring{$K$}{K} is (3)-bad} 
In this subsection, we build towards bounding the density of $\KC$ when $K$ is (3)-bad. Recall that a red component is (3)-bad if it is acyclic, has at most $d$ edges, and has diameter at least $2\ell + 3$.
First, we aim to show that small components are not bad. To that end, we prove the following.

\begin{lemma} \label{lemma:dLargerThanLKPlus1}
    We have that $d > \ell (k+1)$; and if $z \in \{0, \ldots, \ell\}$, we have moreover that  
    $\compSizeBoundWithZ{z} > 2\ell + 1.$
\end{lemma}
\begin{proof}
    We have
    \[\ell (k+1) = \floor{\frac{d-1}{k+1}} (k+1) \leq d - 1.\]
    The second part of the lemma follows immediately since $k \geq 1$.
\end{proof}

This leads us to the useful corollary below.
\begin{cor} \label{cor:smallIsNotBad}
    If $K$ is an acyclic red component with $e(K) \leq 2\ell + 1$, then $K$ is not bad. In particular, small red components are not bad.
\end{cor}

The following technical lemma will be important for several future manipulations we will perform upon $F$.

\begin{lemma} \label{lemma:replaceBadCompWithSmallerBadComps}
    If $K$ is a red acyclic component of $\explSG$ and there is a colouring $f' \in \mathcal F$ with a red pseudoforest $F'$ whose set of components can be obtained from the set of components of $F$ by:
    \begin{itemize}
        \item removing $K$,
        \item possibly adding acyclic components $K_1, ..., K_q$ with $q \in \mathbb N, e(K_i) < e(K)$ and $diam(K_i) \leq diam(K)$ for every $1 \leq i \leq q$, and
        \item possibly adding or removing non-bad components,
    \end{itemize}
    then $K, K_1, ..., K_q$ are not bad and thus $\Delta'(f') \leq_{\text{lex}} \Delta'(f)$.
\end{lemma}
\begin{proof}
    Obviously, the addition and removal of non-bad components does not change $\Delta'$. \\
    First suppose that $K$ is bad. Since $K$ is acyclic, it follows that $K$ is $(b)$-bad for some $b \in \{2, \ldots, 5\}$. Then $\Delta^{(b)}_{e(K)}$ decreases and $\Delta^{(b)}_i$ remains the same for all $i  > e(K)$ when manipulating $F$ as described above. Thus, $\Delta'$ decreases, which is a contradiction.\\
    Next, suppose that $K_i$ is $(b)$-bad, where $i \in \{1, ..., q\}$ and $b \in \{2, ..., 5\}$. But then $K$ is also $(b)$-bad and we again obtain a contradiction as $\Delta'$ decreases when performing the operations described in the lemma.
\end{proof}

In the next two lemmas we will show that the large diameter of $(3)$-bad components prevents them from having small children in our optimal colouring $f$.

\begin{obs} \label{obs:diamToPathLB}
    For any vertex $v$ of a tree $T$, there is a path of length at least $\ceil{\frac{diam(T)}{2}}$ starting at $v$. If $v$ does not lie on a longest path of $T$, then this bound can be increased by one.
\end{obs}
\begin{proof}
    Let $P$ be a path that attains the diameter of $T$. If $v \in V(P)$, then the result is immediate. Otherwise, $v \not \in V(P)$ and let $P' = v, \dots, v'$ be the shortest path from $v$ to $P$. Concatenating $P'$ and the longest path in $P$ starting at $v'$ gives the result. 
\end{proof}

\begin{lemma}
    \label{lemma:boundDiamBySizeOfChild}
    If there is a small child $C$ of a red acyclic component $K$ of $\explSG$ generated by $(x,y)$, then every red path starting at $x$ has length at most $e(C)+1$, and the diameter of $K$ is at most $2e(C) +2$.
\end{lemma}
\begin{proof}
    Suppose towards a contradiction that there is a path inside $K$ of length at least $e(C) + 2$ starting at $x$. Let $P$ be this path, and let $e$ be the edge on this path that is incident with $x$. Exchange $e$ and $(x,y)$. This does not create any cyclic component. Let $K'$ and $K''$ be the resulting new components, where $x \in V(K')$. We have $e(K'') < e(K)$ and $diam(K'') \leq diam(K)$, since $K''$ is a subgraph of $K - e$.\\
    Furthermore, we have that $e(K') \leq e(K) - e(P) + e(C) + 1 < e(K)$ and it is also easy to see that $diam(K') \leq diam(K)$. By Lemma \ref{lemma:replaceBadCompWithSmallerBadComps} we get that $K, K'$ and $K''$ are not bad and thus $\Delta'$ has not increased.\\
    Finally, we can construct a smaller legal order by taking $(R_1, \ldots, R_t)$ up until $K$, and then replacing $K$ with one of $K'$ and $K''$ containing a vertex that witnesses this legal order, and completing the order arbitrarily.
    By this contradiction, it follows that every path of $K$ starting at $x$ has length at most $e(C) + 1$. If $diam(K) \geq 2e(C) + 3$, then by Observation \ref{obs:diamToPathLB} there would be a path of length $e(C) + 2$ starting at $x$, which is a contradiction.
\end{proof}

For red acyclic components of large diameter, Lemma \ref{lemma:boundDiamBySizeOfChild} gives the following.

\begin{cor}	\label{cor:tooLargeDiamMeansNoChildren}
    Let $K$ be an acyclic red component of $\explSG$. If $diam(K) > 2\ell +2$, then $K$ does not have any small children.
\end{cor}

We are now equipped to prove the main result of this subsection.

\begin{lemma}
\label{fracbound3bad}
    If a red component $K$ of $\explSG$ is (3)-bad, then $\densOfKC > \densNDT$.
\end{lemma}
\begin{proof}
    If $K$ is (3)-bad, then by Corollary \ref{cor:tooLargeDiamMeansNoChildren} it does not have any small children. Thus,
    \[
        \densOfKC 
        = \frac{e(K)}{v(K)} 
        \geq \frac{2\ell +3}{2\ell+4} 
        > \frac{\ell + 1}{\ell + 2}
        \geq \densNDT.
    \]
    The last inequality holds due to Lemma \ref{lemma:ellAndDensity}.
\end{proof}

\subsection{Density of \boldmath\texorpdfstring{$\KC$}{K\_C} if \boldmath\texorpdfstring{$K$}{K} is (2)-,(4)- or (5)-bad}

We will not be able to get rid of all the small children of non-small red components that are are not (1)- or (3)-bad like we did in the previous subsections. However, we will bound the number of small children components in the following lemma. A similar lemma can be found in \cite{sndtcPsfs}, but in our case we have to carry out a more careful analysis to ensure that we do not create new bad components when manipulating $f$.

\begin{lemma}
    \label{lemma:atMostKChildren}
    If $K$ is an acyclic red component of $\explSG$, then $K$ has at most $k$ small children. 
\end{lemma}
\begin{proof}
    Suppose that $K$ has $k+1$ small children. As each vertex has only $k$ outgoing blue edges, there exist two distinct vertices $x_1, x_2 \in V(K)$ such that there are two distinct small children $C_1$ and $C_2$ generated by blue arcs $(x_1, y_1)$ and $(x_2, y_2)$, respectively. Let $e_1$ and $e_2$ be the edges on the red path between $x_1$ and $x_2$ that are incident to $x_1$ and $x_2$, respectively. Note that $e_1$ and $e_2$ are not necessarily distinct. For $i \in \{1, 2\}$ let $K_i$ be the component of $K - e_i$ containing $x_i$ and let $K'_i = (V(K_i) \cup V(C_i), \; E(K_i) \cup E(C_i) + x_i y_i)$. Let $L_i = K - V(K_j)$, where $j = 3 - i$. Furthermore, we define \textit{doing exchange i} to mean exchanging $e_i$ and $(x_i, y_i)$.
    In the course of the proof, we will always do only one of the exchanges, not both.
    Note that after doing exchange $i$, the red pseudoforest of the resulting colouring contains the components $K'_i$ and $L_j$ with $e(L_j) < e(K)$. Thus, we call $K'_1, L_2, C_2$ and $K'_2, L_1, C_1$ \textit{new components}. 
    The described components are depicted in Figure \ref{fig:componentsInLeqKChildren}.
    Without loss of generality let $L_2$ contain a vertex that witnesses the legal order (if $K$ is not the root).

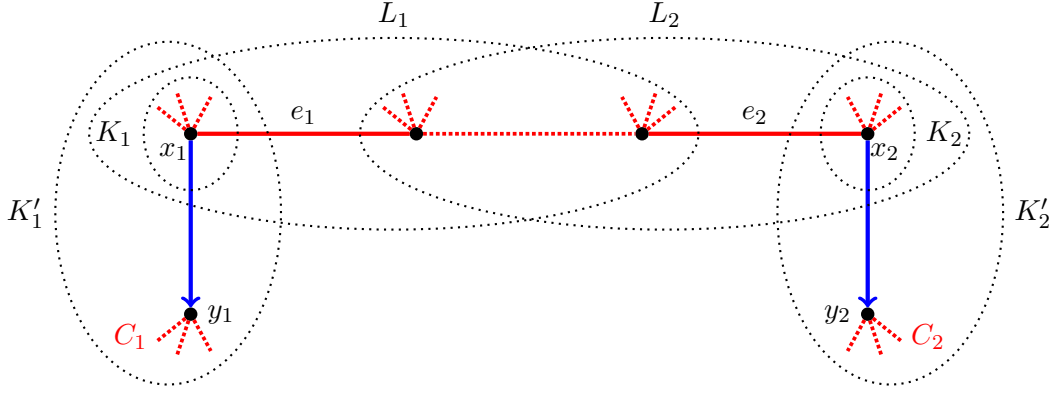
\begin{figure}
\centering
    \begin{tikzpicture}
    \def\sizeOfEdge{3}
    \def\cy{-0.8*\sizeOfEdge}
     \node[blackvertex] at (0,0) (v1) [label = {[label distance=-0.5em]185:$x_{1}$}]{};
     \node[blackvertex] at (1*\sizeOfEdge,0) (v2) {};
     \node[blackvertex] at (2*\sizeOfEdge,0) (v3) {};
     \node[blackvertex] at (3*\sizeOfEdge,0) (v4) [label = {[label distance=-0.5em]-5:$x_{2}$}]{};
     \draw (v1) -- (v2) node [midway, above, fill=white] {$e_{1}$};
     \draw (v3) -- (v4) node [midway, above, fill=white] {$e_{2}$};
     \draw[ultra thick, red] 
        (v1)--(v2)
        (v3)--(v4)
    ;
    \draw[redDotted]
        (v2) -- (v3)
    ;
     \def\transpOneX{-0.5}
     \def\transpOneY{0.4}
     \def\transpTwoX{-0.2}
     \def\transpTwoY{0.6}
     \def\transpThreeX{0.3}
     \def\transpThreeY{0.55}
    \node[transp] at (\transpOneX, \transpOneY) (t11) {};
    \node[transp] at (\transpTwoX, \transpTwoY) (t12) {};
    \node[transp] at (\transpThreeX, \transpThreeY) (t13) {};
    \draw[redDotted]
        (v1) edge (t11)
        (v1) edge (t12)
        (v1) edge (t13)
    ;
    \node[transp] at (1*\sizeOfEdge + \transpOneX, \transpOneY) (t21) {};
    \node[transp] at (1*\sizeOfEdge + \transpTwoX, \transpTwoY) (t22) {};
    \node[transp] at (1*\sizeOfEdge + \transpThreeX, \transpThreeY) (t23) {};
    \draw[redDotted]
        (v2) edge (t21)
        (v2) edge (t22)
        (v2) edge (t23)
    ;
    \node[transp] at (2*\sizeOfEdge - \transpThreeX, \transpThreeY) (t31) {};
    \node[transp] at (2*\sizeOfEdge - \transpTwoX, \transpTwoY) (t32) {};
    \node[transp] at (2*\sizeOfEdge - \transpOneX, \transpOneY) (t33) {};
    \draw[redDotted]
        (v3) edge (t31)
        (v3) edge (t32)
        (v3) edge (t33)
    ;

    \node[transp] at (3*\sizeOfEdge - \transpThreeX, \transpThreeY) (t41) {};
    \node[transp] at (3*\sizeOfEdge - \transpTwoX, \transpTwoY) (t42) {};
    \node[transp] at (3*\sizeOfEdge - \transpOneX, \transpOneY) (t43) {};
    \draw[redDotted]
        (v4) edge (t41)
        (v4) edge (t42)
        (v4) edge (t43)
    ;

     \node[blackvertex] at (0, \cy) (v5)[label = right:$y_{1}$] {};
     \draw[ultra thick, blue, ->] (v1) to (v5);
     \node[blackvertex] at (3*\sizeOfEdge, \cy) (v6) [label = left:$y_{2}$]{};
     \draw[ultra thick, blue, ->] (v4) to (v6);

     \node[transp] at (\transpOneX, \cy -\transpOneY) (ct11) {};
    \node[transp] at (\transpTwoX, \cy -\transpTwoY) (ct12) {};
    \node[transp] at (\transpThreeX, \cy -\transpThreeY) (ct13) {};
    \draw[redDotted]
        (v5) edge (ct11)
        (v5) edge (ct12)
        (v5) edge (ct13)
    ;
    \node[transp] at (3*\sizeOfEdge -\transpThreeX, \cy -\transpThreeY) (c2t1) {};
    \node[transp] at (3*\sizeOfEdge -\transpTwoX, \cy -\transpTwoY) (c2t2) {};
    \node[transp] at (3*\sizeOfEdge -\transpOneX, \cy -\transpOneY) (c2t3) {};
    \draw[redDotted]
        (v6) edge (c2t1)
        (v6) edge (c2t2)
        (v6) edge (c2t3)
    ;
     
     \def\KOffX{0}
     \def\KWidth{1.25 cm}
     \def\KHeight{1.5 cm}
     \def\LOffX{0.9*\sizeOfEdge}
     \def\LWidth{2.7*\sizeOfEdge cm}
     \def\LHeight{1.7*\KHeight}
     \def\CWidth{\KHeight}
     \def\CHeight{\KWidth}
     \def\KPrimeX{-0.3}
     \def\KPrimeY{0.44*\cy}
     \def\KPrimeWidth{1*\sizeOfEdge cm}
     \def\KPrimeHeight{-1.9*\cy cm}
     
    \node[ellipse,draw, thick, dotted, 
        minimum width = \KWidth, 
        minimum height = \KHeight] (e) at (\KOffX,0) [label = left:$K_{1}$] {};
    \node[ellipse,draw,thick,dotted,
    minimum width = \KWidth, 
    minimum height = \KHeight] (e2) at (3*\sizeOfEdge - \KOffX,0) [label = right:$K_{2}$]{};
    \node[ellipse,draw,thick,dotted,
    minimum width = \LWidth, 
    minimum height = \LHeight] (e3) at (\LOffX,0) [label =above:$L_{1}$] {};
    \node[ellipse,draw,thick,dotted,
    minimum width = \LWidth, 
    minimum height = \LHeight] (e4) at (3*\sizeOfEdge - \LOffX,0) [label =above:$L_{2}$] {};
    \node[ellipse,draw,thick,dotted,
    minimum width = \KPrimeWidth, 
    minimum height = \KPrimeHeight] (e5) at (\KPrimeX, \KPrimeY) [label =left:$K_{1}'$] {};
    \node[ellipse,draw,thick,dotted,
    minimum width = \KPrimeWidth, 
    minimum height = \KPrimeHeight] (e6) at (3*\sizeOfEdge - \KPrimeX, \KPrimeY) [label =right:$K_{2}'$] {};
    \node[transp, label={[text=red, label distance=1em]180:$C_1$}] (e7) at (0,\cy - 0.5*\transpTwoY) {};
    \node[transp, label={[text=red, label distance=1em]0:$C_2$}] (e7) at (3*\sizeOfEdge,\cy - 0.5*\transpTwoY) {};
    \end{tikzpicture}
    \caption{A diagram showing the various components in Lemma \ref{lemma:atMostKChildren}. The red dotted lines indicate possibly undrawn vertices and edges. }
    \label{fig:componentsInLeqKChildren}
\end{figure}

    \begin{claim} \label{claim:atMostKChildren:numEdges}
        We have $e(K'_i) < e(K)$ for each $i \in \{1, 2\}$.
    \end{claim}
    \begin{proof}
        Let $i \in \{1,2\}$, $j = 3-i$ and suppose towards a contradiction that we have $e(K'_i) \geq e(K)$. Then we have
        \begin{align*}
            e(K_j) &\leq e(K) - 1 - e(K_i) \\
            &= e(K) - 1 - (e(K'_i) - e(C_i) - 1) \\
            &\leq \ell, \textnormal{\hskip 10mm since $C_i$ is small}.
        \end{align*}
        
        First, suppose that $j = 2$ and $e(K) \leq 2\ell+1$. It follows that $e(K_{1}) \leq \ell$ since $e(K_{2}) \leq \ell$ and $e(K) \geq e(K_{1}) + e(K_{2}) +1$.
        After doing exchange $1$ we have $e(K'_1) \leq 2\ell + 1$ (since $C_1$ is small) and $e(L_2) < e(K)$. Thus, none of the new components are bad by Corollary \ref{cor:smallIsNotBad} and also $K$ is not bad. We see that after the exchange we obtain a smaller legal order by taking the same legal order up to $K$ but replacing $K$ with $L_2$ and then completing the order arbitrarily.\\
        Thus, this case cannot happen and it must be that $j = 1$ or $e(K) > 2\ell + 1$.
        In this case we perform exchange $j$ instead of exchange $i$, then $e(K'_j) \leq 2\ell + 1$, since $e(K_j) \leq \ell$. Note that $K'_j$ is not bad by Corollary \ref{cor:smallIsNotBad}, and so by Lemma \ref{lemma:replaceBadCompWithSmallerBadComps}  we get that neither $K$ nor any of the new components are bad and in particular, $K$ is not the root component.\\
        We can again obtain a smaller legal order by taking the same legal order up to $K$ but replacing $K$ with $K'_2$ if $j=2$ and $e(K) \leq 2\ell + 1$, or with $L_2$ otherwise, and then completing the order arbitrarily.
    \end{proof}
    For each $i \in \{1,2\}$, let $j = 3-i$ and let $r_i$ be the number of edges in a longest path in $K_i$ starting at $x_i$. Note that $r_i \leq e(C_j)$, as otherwise there would be a path of size $e(C_j) +  2$ starting at $x_j$, which is a contradiction to Lemma \ref{lemma:boundDiamBySizeOfChild}. Similarly, $r_i \leq e(C_i) + 1$.
    \begin{claim} \label{claim:atMostKChildren:diam}
        $diam(K'_1), diam(K'_2) \leq 2\ell + 1$.
    \end{claim}
    \begin{proof}
        Suppose $diam(K'_i) \geq 2\ell + 2$.
        As $e(C_i) \leq \ell$ and thus $r_i \leq \ell$, a longest path in $K'_i$ that contains $x_{i}$ has at most $2\ell + 1$ edges and thus, $x_i$ does not lie on a longest path of $K'_i$. But by Observation \ref{obs:diamToPathLB} this gives a path of size $\ell + 2$ starting at $x_i$, which is a contradiction since every path in $K'_i$ starting at $x_i$ is fully contained in either $K_i$ or $C_i + x_iy_i$.
    \end{proof}
    Let us look at the consequences of the two claims regarding whether or not $K$ is bad:
    we defined $K$ to be acyclic, thus it is not (1)-bad and neither are any new components.\\
    Using Claim \ref{claim:atMostKChildren:numEdges} we know that $K$ is not (2)-bad or otherwise $\Delta^{(2)}$ would decrease when doing exchange $1$ or $2$. Thus, none of the new components are (2)-bad.\\
    Analogously, using Claim \ref{claim:atMostKChildren:diam} we could decrease $\Delta'$ if $K$ was (3)-bad or (4)-bad and also none of the new components are (3)-bad or (4)-bad.\\
    If $K$ was (5)-bad, it is again clear by Claim \ref{claim:atMostKChildren:numEdges} that we could decrease $\Delta'$.\\
    If there is an $i \in \{1, 2\}$ such that $K'_i$ is not (5)-bad (and thus not bad, as we excluded all other possibilities), then this either decreases $\Delta'$ or we find a smaller legal order after doing exchange $i$ by taking the same legal order up to $K$ but replacing $K$ with the component $K'_2$, if $i = 2$, or $L_2$, if $i = 1$ and then completing the order arbitrarily.\\
    
    Thus, for the rest of the proof we assume that both $K'_1$ and $K'_2$ are (5)-bad and we aim to show that $K$ is also (5)-bad. This proves the lemma, as in this case we can do either exchange $1$ or $2$ and get a smaller $\Delta'$, by Claim \ref{claim:atMostKChildren:numEdges}, a contradiction.\\
    If a longest path $P$ of $K'_i$ does not contain $x_i$, then either $P \subseteq C_i$ and thus $diam(K'_i) \leq \ell$, or $P \subseteq K_i$, thus $diam(K'_i) = diam(K_i)$ and by Observation \ref{obs:diamToPathLB} we obtain $diam(K'_i) \leq 2r_i$. As $r_{i} \leq e(C_{i}) + 1$, it follows that  $diam(K'_i) \leq r_i + e(C_i) + 1$. On the other hand if there exists a longest path of $K'_i$ containing $x_i$, we also have $diam(K'_i) \leq r_i + e(C_i) + 1$. Thus, in any case we have that $diam(K'_i) \leq r_i + e(C_i) + 1$.

    As both $K'_1$ and $K'_2$ are (5)-bad, there are natural numbers $z_1$ and $z_2$ such that for each $i \in \{1,2\}$, we have $diam(K'_i) \geq 2z_i + 1$ and $e(K'_i) \geq \compSizeBoundWithZ{z_i}$. Thus
    \begin{align*}
        e(K'_i) 
        &\geq \compSizeBoundWithZ{\floor{\frac{diam(K'_i) - 1}{2}}}\\
        &\geq \compSizeBoundWithZ{\floor{\frac{r_i + e(C_i)}{2}}}.
    \end{align*}
    We know that $diam(K) \geq r_1 + r_2 + 1$. Thus, it suffices to prove $e(K) \geq \compSizeBoundWithZ{\floor{\frac{r_1 + r_2}{2}}}$ in order to prove that $K$ is (5)-bad, completing the proof of the lemma as explained above. Observe
    \begin{align*}
        e(K) 
        &\geq e(K_1) + e(K_2) + 1 \\
        &= \sum_{i=1}^2 \big(e(K'_i) - e(C_i) - 1 \big) + 1 \\
        &\geq \sum_{i=1}^2 \Big( d - \floor{\frac{r_i + e(C_i)}{2}} (k-1) - e(C_i) \Big) + 1.
    \end{align*}

From here, let $\alpha_i = 1$ if $r_i \not \equiv e(C_i) \mod 2$ and 0 otherwise. Let $\beta = 1$ if $e(C_1) \not \equiv e(C_2) \mod 2$, and 0 otherwise. Note that $\floor{\frac{r_i + e(C_i)}{2}} = \frac{r_i + e(C_i)-\alpha_i}{2}$. Then 
\begin{align*}
    e(K) &\geq 2d - \left(\frac{e(C_1)+r_1-\alpha_1}{2} + \frac{e(C_2) + r_2 - \alpha_2}{2} \right)(k-1)-(e(C_1) + e(C_2)) + 1 \\
    &= d - \frac{r_1+r_2-(\alpha_1 + \alpha_2)}{2}(k-1)+1+d-\frac{k+1}{2}(e(C_1) + e(C_2)) \\
    &\geq d-\frac{r_1+r_2-(\alpha_1 + \alpha_2)}{2}(k-1)+1+d-\frac{k+1}{2}(2\ell - \beta) \textrm{ since $e(C_i) \leq \ell$ for $i \in \{1,2\}$} \\
    &> d -\frac{r_1 + r_2 - (\alpha_1 + \alpha_2 + \beta)}{2}(k-1) + 1.
\end{align*}
If $1 \in \{\alpha_1, \alpha_2, \beta\}$, then $\floor{\frac{r_1 + r_2}{2}} \geq \frac{r_1 + r_2 - (\alpha_1 + \alpha_2 + \beta)}{2}$. If $ 1 \not \in \{\alpha_1, \alpha_2, \beta\}$, then by definition of $\alpha_1, \alpha_2$, and $\beta$ we have that $r_1 + r_2$ is even, and so again $\floor{\frac{r_1 + r_2}{2}} \geq \frac{r_1 + r_2 - (\alpha_1 + \alpha_2 + \beta)}{2}$. Thus in either case $K$ is (5)-bad since
$$
e(K) \geq d - \floor{\frac{r_1 + r_2}{2}}(k-1)+1.$$
\end{proof}
\begin{lemma}
\label{fracbound2bad}
    If a red component $K$ of $\explSG$ is (2)-bad, then $\densOfKC > \densNDT$.
\end{lemma}
\begin{proof}
    Note that we have $|\mathcal C(K)| \leq k$ by Lemma \ref{lemma:atMostKChildren}. 
    Suppose that $K$ is (2)-bad. Then we have
    \[
        \densOfKC
        = \frac{e(K) + \sum_{C \in \mathcal C(K)} e(C)}{v(K) + \sum_{C \in \mathcal C(K)} v(C)}
        \geq \frac{(d + 1) + k \cdot 0}{(d + 2) + k \cdot 1}
        > \densNDT.
    \]
\end{proof}

In order to show the same for (4)-bad components, we only need the following simple corollary from Lemma \ref{lemma:boundDiamBySizeOfChild}:

\begin{cor}	\label{cor:diamOf2ellPlusOneMeansBigSmallChildrenOfSizeEll}
    Let $K$ be an acyclic red component of $\explSG$. If $diam(K) \geq 2\ell + 1$, then for all small children $C$ of $K$, we have $e(C) = \ell$.
\end{cor}
\begin{lemma}
\label{fracbound4bad}
    If a red component $K$ of $\explSG$ is (4)-bad, then $\densOfKC > \densNDT$.
\end{lemma}
\begin{proof}
    Let $K$ be (4)-bad. Since $\equivModKPlusOne{d}{1}$ we have  $(k+1)\ell = d - 1$. Thus 
        \begin{align*}
            \densOfKC
            &\geq \frac{(2\ell + 2) + k\ell}{(2\ell + 3) + k(\ell + 1)} \\
            &= \frac{(k+1)\ell + \ell + 2}{(k+1)\ell + k + \ell + 3}\\
            &= \frac{d + (\ell + 1)}{d + k + 1 + (\ell + 1)}\\
            &> \densNDT.
        \end{align*}
\end{proof}

\begin{lemma}
\label{fracbound5bad}
    If a red component $K$ of $\explSG$ is (5)-bad, then $\densOfKC > \densNDT$.
\end{lemma}
\begin{proof}
    Suppose $K$ is (5)-bad. Then we have that $e(K) \geq \compSizeBoundWithZ{z}$ and $diam(K) > 2z$ for a $z \in \mathbb N$ with $1 \leq z \leq \ell$. Thus, $e(C) \geq z$ for every small child $C$ of $K$ by Lemma \ref{lemma:boundDiamBySizeOfChild}, which gives us
    \[
        \densOfKC
        \geq \frac{e(K) + kz}{e(K) + 1 + k(z+1)}
        \geq \frac{d - z(k-1) + 1 + kz}{d - z(k-1) + 2 + k(z+1)}
        \geq \frac{d + (z + 1)}{d + k + 1 + (z + 1)}
        > \densNDT.
    \]
\end{proof}

\subsection{Density of \boldmath\texorpdfstring{$\KC$}{K\_C} if \boldmath\texorpdfstring{$K$}{K} is not bad} \label{sec:bad5nonBad}

In this subsection we will show that the density of $K_{C}$ for non-bad components is large if $K$ is not small. If $K$ is small, we will show that $K$ has no small children.
Before this we start with a technical lemma:

\begin{lemma} \label{lemma:diamOfKPlusC}
    Let $K$ be a red acyclic component of $\explSG$ and $C$ be a small child of $K$ generated by $(x, y)$. Then $K' := (V(K) \cup V(C), E(K) \cup E(C) + xy)$ has diameter at most $2e(C) + 2$.
\end{lemma}
\begin{proof}
    If $K'$ contains a path of size $2e(C) + 3$, then by Lemma \ref{lemma:boundDiamBySizeOfChild} this path must contain $xy$ and therefore contain $x$. But again by Lemma \ref{lemma:boundDiamBySizeOfChild}, there are only paths of size at most $e(C) + 1$ starting at $x$ in $K$ and $E(C) + xy$ also has at most $e(C) + 1$ edges.
\end{proof}

\begin{figure}
\begin{center}
    
\begin{tikzpicture}
\node[blackvertex] at (0,0) (v1) {};
\node[blackvertex] at (1.5,0) (v2) {};
\node[blackvertex] at (3,0) (v3) {};
\node[blackvertex] at (4.5,0) (v4) {};
\node[blackvertex] at (0,-1.5) (v5) {};
\node[smallwhite] at (-.5,-1.5) (dummy1) {$x$};
\node[blackvertex] at (0,-3) (v6) {};
\node[smallwhite] at (-.5,-3) (dummy2) {$y$};
\draw[ultra thick, red, -, >=stealth] (v1)--(v2);
\draw[ultra thick,red, -, >= stealth] (v2) --(v3);
\draw[ultra thick, red,-, >= stealth] (v3)--(v4);
\draw[ultra thick, blue, ->, >= stealth] (v1)--(v5);
\draw[ultra thick, blue, ->, >= stealth] (v5)--(v6);
\draw[ultra thick, blue, ->, >= stealth] (v6)--(v3);
\draw[ultra thick, blue, ->, >= stealth] (v3)--(v5);
\draw[ultra thick,blue, ->, bend left = 30, >= stealth] (v2) to (v4);
\draw[ultra thick,blue,->, bend right = 50, >= stealth] (v4) to (v1);

\begin{scope}[xshift = 6cm]
\node[blackvertex] at (0,0) (v1) {};
\node[blackvertex] at (1.5,0) (v2) {};
\node[blackvertex] at (3,0) (v3) {};
\node[blackvertex] at (4.5,0) (v4) {};
\node[blackvertex] at (0,-1.5) (v5) {};
\node[blackvertex] at (0,-3) (v6) {};
\node[smallwhite] at (-.5,-1.5) (dummy1) {$x$};
\node[smallwhite] at (-.5,-3) (dummy2) {$y$};
\draw[ultra thick, blue, ->, >=stealth] (v1)--(v2);
\draw[ultra thick,red, -, >= stealth] (v2) --(v3);
\draw[ultra thick, red,-, >= stealth] (v3)--(v4);
\draw[ultra thick, blue, ->, >= stealth] (v5)--(v1);
\draw[ultra thick, red, -, >= stealth] (v5)--(v6);
\draw[ultra thick, blue, ->, >= stealth] (v6)--(v3);
\draw[ultra thick, blue, ->, >= stealth] (v3)--(v5);
\draw[ultra thick,blue, ->, bend left = 30, >= stealth] (v2) to (v4);
\draw[ultra thick,blue,->, bend right = 50, >= stealth] (v4) to (v1);
\end{scope}
\end{tikzpicture}

\end{center}

       \caption{
            \protect\label{Case2}
            An illustration of Case $2$ in Lemma \ref{lemma:linkKWithCIfBelowZBound}, where $k=1$, $d=2$, $x=w$ and $n=1$. In this case, as the red component had more than two edges, we reduce $\Delta'$. }
\end{figure}
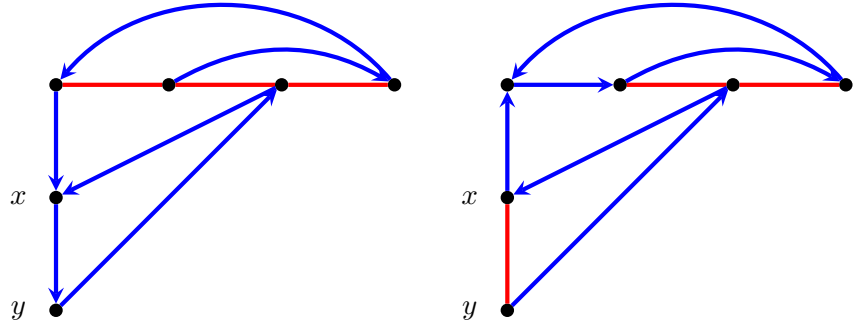

\begin{lemma} \label{lemma:linkKWithCIfBelowZBound}
    Let $K$ be a red component of $\explSG$ that is not bad. Let $C$ be a small child of $K$ generated by $(x, y)$. Furthermore, suppose that $K$ is small or $\nequivModKPlusOne{d}{1}$ or $e(C) < \ell$. Then $e(K) \geq d - e(C)k$.
\end{lemma}
\begin{proof}
Assume to the contrary that $e(K) + e(C) + 1 < \compSizeBoundWithZ{e(C)}$.
Let $L := (V(K) \cup V(C), E(K) \cup E(C) + xy)$.

\begin{claim}\label{claim:Lnotbad}
$L$ is not bad.
\end{claim}
\begin{proof}
Note that $L$ is acyclic and thus not $(1)$-bad. By our assumption we have $e(L) < \compSizeBoundWithZ{e(C)}$ and thus it is not $(2)$-bad. Further, $L$ is not $(5)$-bad. To see this, suppose there exists $z \in \{1, \dots, \ell\}$ such that $e(K) \geq d-z(k-1)+1$. By the assumption, $z > e(C)$. By Lemma \ref{lemma:diamOfKPlusC} it follows that the diameter of $L$ is at most $2e(C)+2$, and thus the diameter of $L$ is at most $2z$ implying that $L$ is not $(5)$-bad. Similarly, by Lemma \ref{lemma:diamOfKPlusC} we obtain that $diam(L) \leq 2e(C) + 2$ and thus it is not (3)-bad.
Finally, we have that $L$ is not (4)-bad, since if $K$ is small we have $e(L) \leq 2\ell + 1$ and if $e(C) < \ell$, then by Lemma \ref{lemma:diamOfKPlusC} we even have $diam(L) \leq 2\ell$.
\end{proof}

Since $K$ is not bad, it is not the root component $R$ and thus there is a vertex $w \in V(K)$ witnessing the legal order.

\textbf{Case 1:} $w \neq x$. \\ 
Let $e$ be the red edge incident to $x$ in $K$ such that $e$ lies on the path from $x$ to $w$ in $K$. Then exchange $e$ and $(x,y)$ and let $K'$ and $K''$ be the new red components containing $x$ and $w$, respectively.
As $K'$ and $K''$ are subgraphs of $L$ and $L$ is not bad by Claim \ref{claim:Lnotbad}, it follows that $K'$ and $K''$ are not bad.
Furthermore, find a smaller legal order by taking the same legal order up to $K$ but replacing $K$ with its proper subgraph $K''$ and completing the order arbitrarily.

\textbf{Case 2:} $w = x$.\\ 
We refer the reader to Figure \ref{Case2} for an illustration. 
As $K$ is not the root component, $K$ has an ancestor.

As $x$ witnesses the legal order, there is a parent component $S_1$ of $K$ that has a blue arc to $x$. If $S_1$ does not have any edges and thus only consists of a single vertex $x_1$, then $x_1$ also witnesses the legal order. In this manner we can find an ancestor of $K$ that contains an edge, since the root component $R$ contains at least one edge. Let $S_1, \ldots, S_n$ be a sequence of red components such that $K$ is a child of $S_1$, $e(S_n) \geq 1$ and for $i \in \{1, \ldots, n-1\}$, $S_i$ is a child of $S_{i+1}$ and $e(S_i) = 0$. There is a blue directed path $P = x_n,\ldots,x_1,x,y$ with $x_{i} \in V(S_{i})$ for all $i \in \{1, \ldots n\}$.
Let $e$ be a red edge incident to $x_{n}$. Now do the following. Colour $(x,y)$ red, remove its orientation and reverse the direction of all remaining arcs in $P$. Colour $e$ blue, and orient $e$ away from $x_{n}$. The resulting coloured mixed graph is in $\mathcal F$, which contains the red and non-bad component $L$. 
By Lemma \ref{lemma:replaceBadCompWithSmallerBadComps} we can conclude that neither $S_{n}$ nor the components of $S_{n}-e$ are bad. Thus, none of the components that have been manipulated by the exchange are bad and in particular, $S_{n}$ is not the root.\\
Finally, we can find a smaller legal order in the colouring of this orientation, as we simply take the same legal order up to $S_{n}$, replace $S_{n}$ with one of the two components of $S_{n}-e$ and complete the remaining order arbitrarily.
\end{proof}

\begin{cor}
    \label{smallchildrencorollary}
    If $K$ is a small red component of $\explSG$, then $K$ does not have any small red children. Furthermore, every small red component of $\explSG$ has a parent component which is not small.
\end{cor}
\begin{proof}
    Suppose $K$ is small and has a small child $C$. Then $e(K) \leq \ell$ by definition of small, and by Lemma \ref{lemma:linkKWithCIfBelowZBound} $e(K) \geq d-e(C)k$. Moreover since $C$ is small, $e(K) \geq d-\ell k$. Thus $d-\ell k \leq \ell$, or $d \leq \ell(k+1)$. This contradicts the definition of $\ell$. 
    For the second part of the corollary remember that the root component is bad and therefore not small due to Corollary \ref{cor:smallIsNotBad}.
\end{proof}

\begin{lemma}
\label{fracboundnotbad}
    If a red component $K \in \mathcal K$ of $\explSG$ is not bad, then $\densOfKC \geq \densNDT$.
\end{lemma}
\begin{proof}
    First, let $\equivModKPlusOne{d}{1}$ and thus $(k+1)\ell = d - 1$. Suppose all small children of $K$ have $\ell$ edges and recall that $K$ has at most $k$ small children by Lemma \ref{lemma:atMostKChildren}. Then
    \[
        \densOfKC 
        \geq \frac{e(K) + k\ell}{e(K) + 1 + k(\ell + 1)}
        \geq \frac{\ell + 1 + k\ell}{\ell + 2 + k(\ell + 1)}
        = \frac{(k+1)\ell + 1}{(k+1)\ell + k + 2}
        = \densNDT.
    \]
    
    Otherwise, we can apply Lemma \ref{lemma:linkKWithCIfBelowZBound} and obtain $e(K) \geq d - ke(C)$, where $C$ is a small child of $K$ with the smallest number of edges. In this case we have
    \[
        \densOfKC 
        \geq \frac{d - ke(C) + ke(C)}{d - ke(C) + 1 + k(e(C) + 1)}
        \geq \densNDT.
    \]
\end{proof}

Theorem \ref{thm:upperBound} now follows from Lemma \ref{lemma:howToGetMADContradiction}. 

\def\numSimpleComponents{p}

\section{Lower Bound of the Overall Diameter}

We recall the statement of Theorem \ref{thm:lowerBound}.

\restlowerBound*

In this section, we will construct a graph $G$ parameterized by $\delta$ and $\numSimpleComponents$ that satisfies Theorem \ref{thm:lowerBound} for fixed $k, D, \epsilon, \alpha$. We do this in phases. First we define a specific tree with a particular orientation and edge-colouring. Given a tree, we will call the set of vertices that are not leaves \textit{inner vertices} of the tree.

Let $T$ be a tree with root $r_T$ and odd depth $\delta \geq \frac{2(\ell + 1)}{\epsilon} - 1$ formed from a 1-coloured directed path of length $\delta$ starting at $r_T$ by adding another $k-1$ outgoing edges (to new vertices) with colours $2, \ldots, k$ to each even-depth vertex of this path.

From $T$, we construct a tree $C$ in the following manner. For each $t \in V(T)$, we add a path $P_{t}$ which is disjoint from $T$ except for at $t$, and $P_{t}$ has length $l_t$ where

\[
    l_t = 
    \begin{cases*}
        2\ell + 1 + \alpha & \text{if $t = r_T$,}\\
        \ell + \alpha & \text{if $t \neq r_T$ and $depth_{T}(t)$ is even,}\\
        \ell & \text{if $depth_{T}(t)$ is odd.}\\
    \end{cases*}
\]
If $t \neq r_T$, then $t$ is an endpoint of $P_t$; and $r_T$ is in the middle of $P_{r_T}$ (i.e. there are two edge-disjoint subpaths of $P_{r_T}$ starting at $r_T$ with length at least $\ell$).
Further, for all $t \in V(T)$, we colour all edges of $P_{t}$ with $k+1$, and we orient the edges towards $t$. Let us call $C$ a \textit{colourful tree}. We say that $r_C := r_T$ is the root of $C$ as well.

Finally, we obtain our desired graph $G$ by taking $p$ pairwise disjoint copies of $C$, where $p$ is any integer such that $p \geq kD + k^2 + 1$, a set of new vertices $S := \{s_{1},\ldots,s_{k}\}$ and for every colour $i \in \{1, \ldots k\}$ and every vertex $x$ belonging to a colourful tree we add an edge $(x, s_i)$, if $x$ does not have an outgoing edge coloured $i$. For each copy of $C$, we let $T_{C}$ denote that copy of $T$ contained in $C$.

We will refer to the colours $\{1,\ldots,k\}$ as blue, and the colour $k+1$ as red. Note that unlike for the upper bound, where it is useful to think of all of the different blue colours as the same, for this construction each blue monochromatic component should be thought of as having a specific blue colour $i$. We call the established colouring and orientation of $G$ the \textit{example colouring}. A colourful tree having the colours of the example colouring is depicted in Figure \ref{fig:exampleColouring}. Note the example colouring is a colouring where with the tools developed in Section 2, we would not be able to reduce the diameter further, and the red pseudoforest contains a path of length $2\ell + 1 + \alpha$. The example colouring will be especially useful in Subsection \ref{subsec:boundingFracArb} as it allows us to concisely refer to different structures within $G$.

Our first point of order is to show that in any decomposition of $G$ into $k+1$ pseudoforests, the red pseudoforest has large diameter. We do this in Subsection \ref{subsec:boundingDiam}. In Subsection \ref{subsec:boundingFracArb}, we lower-bound the fractional arboricity of $G$, which completes the proof of Theorem \ref{thm:lowerBound}.

\subsection{Bounding the Diameter}\label{subsec:boundingDiam}

\begin{thm} \label{thm:lowerBoundDiam}
    In any decomposition of $G$ into $k$ blue pseudoforests and one red pseudoforest where every vertex has at most $D$ incident red edges, there is a red component that has diameter at least $2\ell + 1 + \alpha$.
\end{thm}

We assume there is a colouring $f$ that contradicts Theorem \ref{thm:lowerBoundDiam}. In the following lemma and corollary we will easily force colours on many edges of at least one colourful tree. The result is depicted in Figure \ref{fig:forcedColoursOnEdges}. The rest of the subsection will show that we cannot find colours for the remaining black edges in the figure without creating red paths which are too long or creating red components with too many cycles.\\

\begin{figure}
\centering
\begin{minipage}[t]{.5\textwidth}
  \centering
    
\begin{tikzpicture}
\node[blueDot]  (v00) at (-2*\edgeSizeLB, 0) {};
\node[blueDot]  (v01) at (-1*\edgeSizeLB, 0) {};
\node[blackDot] (v02) at ( 0*\edgeSizeLB, 0) [label = {90:$r_C$}] {};
\node[blueDot]  (v03) at ( 1*\edgeSizeLB, 0) {};
\node[blueDot]  (v04) at ( 2*\edgeSizeLB, 0) {};
\node[blueDot]  (v10) at (-1*\depthOneX - 1*\edgeSizeLB, \depthOneY) {};
\node[blueDot] (v11) at (-1*\depthOneX - 0*\edgeSizeLB, \depthOneY) {};
\node[cyanDot] (v12) at ( 0*\depthOneX + 0*\edgeSizeLB, \depthOneY) {};
\node[blueDot]  (v13) at ( 0*\depthOneX + 1*\edgeSizeLB, \depthOneY) {};
\node[blueDot] (v14) at ( 1*\depthOneX + 0*\edgeSizeLB, \depthOneY) {};
\node[blueDot]  (v15) at ( 1*\depthOneX + 1*\edgeSizeLB, \depthOneY) {};
\node[blackDot]  (v20) at ( 0*\edgeSizeLB, \depthTwoY) {};
\node[blueDot]  (v21) at ( 1*\edgeSizeLB, \depthTwoY) {};
\node[blueDot] (v22) at ( 2*\edgeSizeLB, \depthTwoY) {};

\node[blueDot] (v30) at (-\depthThreeX - 1*\edgeSizeLB, \depthThreeY) {};
\node[blueDot] (v31) at (-\depthThreeX - 0*\edgeSizeLB, \depthThreeY) {};

\node[blueDot] (v32) at (0*\edgeSizeLB, \depthThreeY) {};
\node[blueDot] (v33) at (1*\edgeSizeLB, \depthThreeY) {};

\node[blueDot] (v34) at ( \depthThreeX + 0*\edgeSizeLB, \depthThreeY) {};
\node[blueDot] (v35) at ( \depthThreeX + 1*\edgeSizeLB, \depthThreeY) {};

\draw[redEdge]
    (v00) edge (v01)
    (v01) edge (v02)
    (v02) edge (v03)
    (v03) edge (v04)
    (v10) edge (v11)
    (v12) edge (v13)
    (v14) edge (v15)
    (v20) edge (v21)
    (v21) edge (v22)
    (v30) edge (v31)
    (v32) edge (v33)
    (v34) edge (v35)
    ;
\draw[blueEdge]
    (v02) edge (v11)
    (v02) edge (v12)
    (v02) edge (v14)
    (v12) edge (v20)
    (v20) edge (v31)
    (v20) edge (v32)
    (v20) edge (v34)
    ;
\end{tikzpicture}
  \caption{The example colouring of a colourful tree if $k = 3, \ell = 1, \alpha = 1, \delta = 3$. Blue (light blue) vertices have an edge to $S$ in every colour of $\{1, \ldots, k\}$ (of $\{2, \ldots k\})$. }
  \label{fig:exampleColouring}
\end{minipage}%
\begin{minipage}[t]{.5\textwidth}
  \centering
    \begin{tikzpicture}
\node[blueDot]  (v00) at (-2*\edgeSizeLB, 0) {};
\node[blueDot]  (v01) at (-1*\edgeSizeLB, 0) {};
\node[blackDot] (v02) at ( 0*\edgeSizeLB, 0) [label = {90:$r_C$}] {};
\node[blueDot]  (v03) at ( 1*\edgeSizeLB, 0) {};
\node[blueDot]  (v04) at ( 2*\edgeSizeLB, 0) {};
\node[blueDot]  (v10) at (-1*\depthOneX - 1*\edgeSizeLB, \depthOneY) {};
\node[blueDot] (v11) at (-1*\depthOneX - 0*\edgeSizeLB, \depthOneY) {};
\node[cyanDot] (v12) at ( 0*\depthOneX + 0*\edgeSizeLB, \depthOneY) {};
\node[blueDot]  (v13) at ( 0*\depthOneX + 1*\edgeSizeLB, \depthOneY) {};
\node[blueDot] (v14) at ( 1*\depthOneX + 0*\edgeSizeLB, \depthOneY) {};
\node[blueDot]  (v15) at ( 1*\depthOneX + 1*\edgeSizeLB, \depthOneY) {};
\node[blackDot]  (v20) at ( 0*\edgeSizeLB, \depthTwoY) {};
\node[blueDot]  (v21) at ( 1*\edgeSizeLB, \depthTwoY) {};
\node[blueDot] (v22) at ( 2*\edgeSizeLB, \depthTwoY) {};

\node[blueDot] (v30) at (-\depthThreeX - 1*\edgeSizeLB, \depthThreeY) {};
\node[blueDot] (v31) at (-\depthThreeX - 0*\edgeSizeLB, \depthThreeY) {};

\node[blueDot] (v32) at (0*\edgeSizeLB, \depthThreeY) {};
\node[blueDot] (v33) at (1*\edgeSizeLB, \depthThreeY) {};

\node[blueDot] (v34) at ( \depthThreeX + 0*\edgeSizeLB, \depthThreeY) {};
\node[blueDot] (v35) at ( \depthThreeX + 1*\edgeSizeLB, \depthThreeY) {};

\draw[redEdge]
    (v00) edge (v01)
    (v03) edge (v04)
    (v10) edge (v11)
    (v14) edge (v15)
    (v21) edge (v22)
    (v30) edge (v31)
    (v32) edge (v33)
    (v34) edge (v35)
    ;
\draw[blackEdge]
    (v01) edge (v02)
    (v02) edge (v03)
    (v02) edge (v11)
    (v02) edge (v12)
    (v02) edge (v14)
    (v12) edge (v13)
    (v12) edge (v20)
    (v20) edge (v21)
    (v20) edge (v31)
    (v20) edge (v32)
    (v20) edge (v34)
    ;
\end{tikzpicture}
  \caption{$C$ after enforcing the edge-colouring of Corollary \ref{cor:forcedColourOfEdges}. Blue (light blue) vertices have $k$ ($k-1$) blue edges to $S$.}
  \label{fig:forcedColoursOnEdges}
\end{minipage}
\end{figure}

In what follows, an \textit{$S$-$C$-$S$-path} is a path with endpoints in $S$ and whose inner vertices are all from one colourful tree $C$.

\begin{lemma} \label{lemma:colourfulEdgesToS}
    There is a colourful tree $C$ in $G$ such that every edge of $E(C, S)$ is coloured blue in $f$. Furthermore, there is no monochromatic $S$-$C$-$S$-path.
\end{lemma}
\begin{proof}
    There are at most $kD$ red edges incident to $S$. For any colour $b \in \{1, \ldots, k\}$ there can be at most $k$ $S$-$C$-$S$ paths having colour $b$ in $f$ as otherwise we would have monochromatic components containing two cycles.
    As $\numSimpleComponents > kD + k^2$ there is at least one colourful tree $C$ which satisfies the lemma.
\end{proof}

For the rest of the subsection we let $C$ be a colourful tree satisfying Lemma \ref{lemma:colourfulEdgesToS}. The following corollary follows easily from Lemma \ref{lemma:colourfulEdgesToS}.

\begin{cor}	\label{cor:forcedColourOfEdges}
    In $f$ any vertex of $C$ that is not an inner vertex of $T_C$ has a $b$-coloured edge to $S$ for every $b \in \{1, \ldots, k\}$
    and any vertex of $T_C$ with odd depth has $k-1$ incident edges to $S$ in pairwise distinct colours of $\{1, \ldots, k\}$.
    Any edge of $E(C)$ that is not incident with an inner vertex of $T_C$ is coloured red in $f$. Furthermore, every red component containing at least one vertex in $C$ is acyclic in $f$.
\end{cor}

Given a colour $b \in \{1, \dots, k\}$, we say a vertex $v \in V(C)$ has a \textit{$b$-coloured $S$-path} if there is a path of colour $b$ in $f$ that goes from $v$ to a vertex of $S$ and for any inner vertex $w$ of this path, $w \in V(C)$ and $depth_C(w) \geq depth_C(v)$.
Further, we say a vertex \textit{$t$ is the end of a low $i$-path} if $t$ is an endpoint of a red path $P$ in $f$ with $e(P) \geq i$ and $depth_C(v) \geq depth_C(t)$ for every $v \in V(P)$.

For $i \in \{\ell, \ell + 1, 2\ell + 1 + \alpha\}$, let $V_i$ be the set of vertices of $V(T_C)$ contained in red components with exactly $i$ edges in the example colouring.

\begin{lemma} \label{lemma:inductionOverNonRootVertices}
    In $f$, for each $i \in \{\ell, \ell+1\}$ every $t \in V_i$ is the end of a low $i$-path and has a $b$-coloured $S$-path for every $b \in \{1, \ldots, k\}$.
\end{lemma}
\begin{proof}
    The lemma is clear for any leaf of $T_C$ due to Corollary \ref{cor:forcedColourOfEdges}. Next, suppose that $t \in V_{\ell + \alpha}$ 
    has even depth, is not the root of $T_C$ and that the lemma is true for all of its
    children $c_1, \ldots, c_k$ in $T_C$. Let $u$ be the vertex such that there is a red arc $(u, t)$ in the example colouring. Note that $u$ is also a child of $t$ in $C$. By Corollary \ref{cor:forcedColourOfEdges}, $u$ is the end of a low $(\ell + \alpha - 1)$-path in $f$ and it also has a $b$-coloured $S$-path (of length 1) for any $b \in \{1, \ldots, k\}$.\\
    No two edges of $t c_1, \ldots, t c_k, t u$ have the same colour  $j \in \{1, \ldots, k\}$ in $f$ or there would be a monochromatic $S$-$C$-$S$-path, contradicting Lemma \ref{lemma:colourfulEdgesToS}. Furthermore, no two of these edges are red or $f$ would have a red path with at least $(\ell + \alpha - 1) + 2 + \ell$ edges, a contradiction.
    Thus the $k+1$ edges have pairwise distinct colours, from which the lemma follows.

    Now, let $t \in V_\ell$ have odd depth, not be a leaf of $T_C$ and again suppose that the lemma is true for the only child $t'$ of $t$ in $T_C$. Let $u$ be the other child of $t$ in $C$, which is the neighbour of $t$ in $P_t$. Note that $t'$ has a low $(\ell + \alpha)$-path, and by Corollary \ref{cor:forcedColourOfEdges} $u$ has a low $(\ell-1)$-path.
    Since $t$ has $k-1$ $S$-paths of length $1$, at most one of the edges to the children of $t$ can be blue.
    If both edges were red, then there would be a red path of length $(\ell+\alpha) + 2 + (\ell-1)$, which is contrary to our assumptions. The lemma follows.
\end{proof}

The following corollary shows that $r_C$ is contained in a red component with too high diameter, contradicting our initial assumption and thus completing the proof of Theorem \ref{thm:lowerBoundDiam}.

\begin{cor}
    The diameter of the red component of $r_C$ in $f$ is at least $2\ell+1 + \alpha$.
\end{cor}
\begin{proof}
    By Lemma \ref{lemma:inductionOverNonRootVertices} we know that the $k$ children of $r_C$ in $T_C$ are the ends of low $\ell$-paths and that they each have a $b$-coloured $S$-path for all $1 \leq b \leq k$. Furthermore, by Corollary \ref{cor:forcedColourOfEdges} the two neighbours of $r_C$ in $P_{r_C}$ are also ends of low $\ell$- and $(\ell + \alpha - 1)$-paths, respectively, and they have $b$-coloured $S$-paths of length $1$ for every $1 \leq b \leq k$.\\
    We again conclude that there cannot be two edges from $r_C$ to one of its $k+2$ children in $C$ that have the same colour $b$ with $1 \leq b \leq k$ as otherwise we have a monochromatic $S$-$C$-$S$-path, contradicting Lemma \ref{lemma:colourfulEdgesToS}. Thus, two of the $k+2$ edges from $r_T$ to its children are red and thus \textemdash{} together with the low $\ell$- (and perhaps $(\ell + \alpha -1)$-) paths described above\textemdash{} form a red path with at least $(\ell + \alpha - 1) + \ell + 2 = 2\ell + 1 + \alpha$ edges.
\end{proof}

\subsection{Bounding the Fractional Arboricity}\label{subsec:boundingFracArb}

The last point of order is to bound the fractional arboricity of $G$. As mentioned in the introduction, this will also give us a bound on its maximum average degree.
We aim to show that the entire graph is the densest subgraph of $G$. We have chosen $\delta$ large enough such that the red components of a colourful tree in the example colouring that do not contain a root of a colourful tree, which have a density of roughly $\lowerMadBoundWithAlpha$, compensate for the largest red component (which contains the root).

\begin{thm} \label{thm:lowerBoundFracArb}
    The fractional arboricity of $G$ is less than $k + \lowerMadBoundWithAlpha+ \epsilon$.
\end{thm}

We assume towards a contradiction that there is a subgraph $H \subseteq G$ with $v(H) \geq 2$ and $\frac{e(H)}{v(H) - 1} \geq k + \lowerMadBoundWithAlpha + \epsilon$  for the rest of the subsection. Let $\mathcal X$ be the set of all subsets $X \subseteq V(G)$ with  $S \subsetneq X$ and 
\[\frac{e(X)}{|X \setminus S|} > k + \beta,\]
where
\[
    \beta = \frac{ \ell(k+1) + \alpha + \frac{2(\ell+1)}{\delta + 1} }
    { (\ell+1)(k+1) + \alpha + \frac{2(\ell+1)}{\delta + 1} }
\]

\begin{lemma}
\label{contradictionforbound}
    $\frac{e(G)}{|V(G) \setminus S|} = k + \beta$ and thus $V(G) \not\in \mathcal X$.
\end{lemma}
\begin{proof}
    In the example colouring every colourful tree $C$ has $\frac{\delta + 1}{2}$ red components $P_t$ containing $\ell$ edges (where $depth_{T_C}(t)$ is odd) and $C$ has $\frac{\delta - 1}{2}$ red components $P_t$ containing $\ell + \alpha$ edges (where $depth_{T_C}(t)$ is even and non-zero).
    Furthermore, note that every vertex of $V(G) \setminus S$ has exactly $k$ outgoing blue edges in the example colouring and there are no outgoing blue edges from $S$. It follows that
    \[
        \frac{e(G)}{|V(G) \setminus S|}
        = k + \frac{
            \numSimpleComponents \big(
                2\ell + \alpha + 1 + \frac{\delta - 1}{2} (\ell + \alpha) + k  \frac{\delta + 1}{2} \ell 
            \big)
        }{
            \numSimpleComponents \big(
                2\ell + \alpha + 2 + \frac{\delta - 1}{2} (\ell + \alpha + 1) + k  \frac{\delta + 1}{2} (\ell + 1) 
            \big)
        }
        = k + \beta.
    \]
    
\end{proof}

\begin{lemma}
     \[
        \lowerMadBoundWithAlpha 
        < \beta 
        < \lowerMadBoundWithAlpha + \epsilon.
     \]
\end{lemma}
\begin{proof}
    The first inequality is easy to see. For the second inequality we use the lower bound of $\delta \geq \frac{2(\ell + 1)}{\epsilon} - 1$ and the fact that for any numbers $a \geq 0, b \geq 1$, and $\gamma > 0$, it holds that $\frac{a + \gamma}{b + \gamma} < \frac{a}{b} + \gamma$.\\
\end{proof}

\begin{cor} \label{cor:connectionWeirdDensityToFracArb}
    If Theorem \ref{thm:lowerBoundFracArb} is false, then $\mathcal X$ is not empty.
\end{cor}
\begin{proof}
    Let $H$ be an induced subgraph of $G$ with $v(H) \geq 2$ and $\frac{e(H)}{v(H) - 1} \geq k + \lowerMadBoundWithAlpha + \epsilon$. Note that $S \subseteq V(H)$, as otherwise $H$ decomposes into at most $k-1$ star forests, where the centers lie in $S$, and one forest whose edges are from colourful trees, contradicting Theorem  \ref{nashthm} since $\frac{e(H)}{v(H)-1} > k$ by assumption. We have
    \[
        k + \beta
        < k + \lowerMadBoundWithAlpha + \epsilon 
        \leq \frac{e(H)}{v(H) - 1} 
        \leq \frac{e(H)}{|V(H) \setminus S|}.
    \]
\end{proof}

Our goal is to show that $V(G) \in \mathcal X$, which contradicts Lemma \ref{contradictionforbound}. For this we will show that if there is an $X \in \mathcal X$, then we can manipulate $X$ in multiple steps such that after every step the new set of vertices is still in $\mathcal X$, and the resulting set is $V(G)$. The technical lemma that we will use implicitly in the remaining lemmas is the following, which follows from the definition of $\mathcal{X}$:

\begin{lemma}
Let $X \in \mathcal X$ and $Z \subseteq V(G) \setminus S$.\\
If $Z \subseteq X$ and 
\[
    \frac{ e(Z) + e(X \setminus Z, Z) }{|Z|} \leq k + \beta,
\]
then $X \setminus Z \in \mathcal X$ and in particular, $X \setminus Z \neq S$.
If $X \cap Z = \varnothing$ and 
\[
    \frac{e(Z) + e(X, Z)}{|Z|} \geq k + \beta,
\]
then $X \cup Z \in \mathcal X$.
\end{lemma}

Let $C$ be a colourful tree and $t \in V(T_C)$ for the rest of the subsection.

\begin{lemma}   \label{lemma:wholePInX}
    If $X \in \mathcal X$ and $t \in X$, then $X \cup V(P_t) \in \mathcal X$.
\end{lemma}
\begin{proof}
    If we add all vertices of $V(P_t) \setminus X$ to $X$, then for each vertex $v \in V(P_t) \setminus X$ in the induced subgraph we add at least $k$ edges from $v$ to $S$ and at least one other edge, namely the red outgoing edge of $v$ in the example colouring. Thus, 
    \[
        e(V(P_t) \setminus X) + e(V(P_t) \setminus X, X) \geq (k+1) \; |V(P_t) \setminus X|.
    \]
\end{proof}

\begin{lemma}
    If $X \in \mathcal X$ and $t \not\in X$, then $X \setminus V(P_t) \in \mathcal X$.
\end{lemma}
\begin{proof}
    Let $P' = X \cap V(P_t) \neq \varnothing$. When subtracting $P'$ from $X$ the induced subgraph will lose $k|P'|$ edges of $E(P', S)$ and another $e(P') \leq |P'| - 1$ edges. If $t \neq r_C$ or $G[P']$ has exactly one component, then we have
    \[
        \frac{e(P')+e(P',S)}{|P'|} = \frac{k|P'| + e(P')}{|P'|} 
        \leq k + \frac{|P'| - 1}{|P'|} 
        \leq k + \frac{\ell}{\ell + 1}
        < k + \beta.
    \]
    If $t = r_C$ and $G[P']$ has two components, we have that $e(P') \leq |P'| - 2$ and $|P'| \leq 2\ell + 1 + \alpha$. In this case we have
    \[
        \frac{e(P')+e(P',S)}{|P'|} =\frac{k|P'| + e(P')}{|P'|} 
        \leq k + \frac{|P'| - 2}{|P'|} 
        \leq k + \frac{\ell + (\ell + \alpha - 1)}{\ell + 1 + (\ell + \alpha)} 
        \leq k + \frac{\ell}{\ell + 1}
        < k + \beta.
    \]
\end{proof}

\begin{lemma} \label{lemma:ellHasDegkPlusTwo}
    If $X \in \mathcal X$, $t \in X \cap V_\ell$ and a neighbour of $t$ in $T_C$ is not in $X$, then $X \setminus V(P_t) \in \mathcal X$.
\end{lemma}
\begin{proof}
    Let $X' = X \cup V(P_t)$. We have $X' \in \mathcal X$ due to Lemma \ref{lemma:wholePInX}. When removing the $\ell + 1$ vertices of $P_t$ from $G[X']$, we remove $k\ell$ edges of $E(V(P_t) - t, S)$, at most $k+1$ edges that were incident to $t$ and all $\ell -1$ edges of $E(P_t -t)$. We have
    \[
        \frac{k\ell + (k+1) + (\ell - 1)}{\ell + 1}
        = \frac{k(\ell + 1) + \ell}{\ell + 1} 
        = k + \frac{\ell}{\ell + 1}
        < k +\beta
    \]
    and thus $X \setminus V(P_t) = X' \setminus V(P_t) \in \mathcal X$.

\end{proof}

By repeated application of the last three lemmas we get the following corollary:

\begin{cor}\label{cor:Xexists}
    If $\mathcal X \neq \varnothing$, then there is an $X \in \mathcal X$ such that for every $t \in V(T_C) \cap X$ we have $V(P_t) \subseteq X$ and if additionally $t \in V_\ell$, then $deg_{X}(t) = k+2$.
\end{cor}

We also want to prove the degree property of Corollary \ref{cor:Xexists} for every $t \in V(T_C)$. We will tackle this in the following two lemmas.

\begin{lemma} \label{lemma:noPerfectKAryTrees}
    Let $X \in \mathcal X$ such that for every $t \in V(T_C) \cap X$ we have $V(P_t) \subseteq X$ and if additionally $t \in V_\ell$, then $deg_{X}(t) = k+2$.
    Furthermore, suppose there is a vertex $r' \in V(T_C)$ such that all vertices of the tree $T'$ containing the vertices $t \in V(T_C)$ with $depth(t) \geq depth(r')$ are also in $X$, but the parent of $r'$ is not in $X$. Let $T'' = T' \cup \bigcup_{t \in V(T')} P_{t}$. Then $X \setminus V(T'') \in \mathcal X$.
\end{lemma}
\begin{proof}
    We have that $r' \in V_{\ell+1}$, as $r'$ is the root of $T'$ and thus has at most $k+1$ neighbours in $X$, and all vertices in $V_{\ell} \cap X$ have $k+2$ neighbours in $X$. It follows that $\alpha =1$ since $V_{\ell+1} \neq \varnothing$.
    Let $T''$ contain $x$ vertices of $V_{\ell+1}$ and thus $kx$ vertices of $V_\ell$. We have
   \[
      \frac{e(T'') + e(X \setminus V(T''), V(T''))}{v(T'')}
       = k + \frac{ x(\ell+1) + kx\ell }{ x(\ell+2) + kx(\ell+1) }
        = k + \lowerMadBoundWithAlpha
        < k + \beta.
   \]

   Thus, $X \setminus V(T'') \in \mathcal X$.
\end{proof}

\begin{lemma} \label{lemma:ellPlusOneHasKChildren}
    If $\mathcal X \neq \varnothing$, then there is an $X \in \mathcal X$ where every $t \in V_\ell \cap X$ has $deg_{X}(t) = k+2$ and for every $t \in X \cap (V_{\ell+1} \cup V_{2\ell + 1 + \alpha})$ the child of $t$ in $T_C$ is also in $X$.
\end{lemma}
\begin{proof}
    Note that by Corollary \ref{cor:Xexists}, there exists an $X \in \mathcal{X}$ where every $t \in V_\ell \cap X$ has $\deg_X(t) = k+2$. From all such $X$, we choose one set $X$ where
    the number of vertices $t \in X \cap (V_{\ell+1} \cup V_{2\ell + 1 + \alpha})$, where the child of $t$ in $T_C$ is not in $X$, is minimal.
    Suppose towards a contradiction that this minimum value is greater than zero and let $t$ be such a vertex with maximal depth and its child being $t' \not\in X$.\\
    Let $T'$ be the subtree of $C$ induced by $t'$ and all of its descendants. By Lemma \ref{lemma:noPerfectKAryTrees} we can choose $X$ such that no vertex of $T'$ is in $X$. Our desired contradiction will be $X' := X \cup V(T') \cup \bigcup_{t\in V(T')} V( P_{t} ) \in \mathcal X$, since this set contains $t'$. Let $x := |V_{\ell+1} \cap V(T')|$ and thus $|V_\ell \cap V(T')| = kx+ |\{t'\}|$. Since $T' \neq T_C$, we have $ x < \frac{\delta - 1}{2}$ and thus

    \[
        \frac{e(X')}{|V(X') \setminus S|}
        = k + \frac
            {
                (\ell + |\{t t'\} | ) + x(\ell + \alpha) + kx\ell
            }{
                (\ell + 1) + x(\ell + \alpha + 1) + kx(\ell + 1)
            }
        = k + \frac
            {
                \ell(k+1) + \alpha + \frac{\ell+1}{x}
            }{
                (\ell+1)(k+1) + \alpha + \frac{\ell+1}{x}
            }\\
        > k + \beta.
    \]
\end{proof}

\begin{cor}
    $\mathcal X$ is empty. Thus, Theorem \ref{thm:lowerBoundFracArb} is true.
\end{cor}
\begin{proof}
    Let $X \in \mathcal X$ such that it satisfies Lemma \ref{lemma:ellPlusOneHasKChildren}. We can choose $X$ such that its induced subgraph does not have any subtrees $T'$ with root $r'$ as described in Lemma \ref{lemma:noPerfectKAryTrees}. But then we have either $X = S$ or $X = V(G)$, both of which give a contradiction: the former to the definition of $\mathcal{X}$, and the latter to Lemma \ref{contradictionforbound}.
\end{proof}

\section{Lower Bound of the Diameter of Big Components}

In this section we will prove Theorem \ref{thm:z:lowerBound}, restated below for convenience.

\restzlowerBound*

To do this, we use a very similar construction to the construction given in the last section. The structure of this section is also identical to the last section and for some proofs we will just refer to the previous section. First we give the construction. We will again describe the \textit{example colouring of $G$}:

For every colourful tree $C$ the tree $T_C$ has odd depth $\delta > \frac{2(z+1)}{\epsilon} - 1$,  and is comprised of a directed blue path of colour $1$ and length $\delta - 1$ starting at the root $r_C$ of $T_C$ and $C$, where additionally, every vertex with even depth has $k-1$ other outgoing edges with colours $2, \ldots, k$ in $T_C$. These $k-1$ children are leaves in $T_C$. We obtain $C$ from $T_C$ by adding the fewest number of new vertices and red edges such that no two vertices of $V(T_C)$ are in the same red component and for every vertex $t \in V(T_C)$ we have:
\begin{itemize}
    \item If $depth(t)$ is odd, then there is an induced red path of length $z$ ending in $t$. 
    \item If $t = r_C$, then $t$ is in a red acyclic component $K$ with $e(K) = d - z(k-1) + 1$ and $diam(K) = 2(z+1)$, $t$ is exactly in the middle of a red path $Q_t \subset K$ of length $2(z+1)$ and $K - t$ consists of components with at most $z$ edges.
    \item If $depth(t)$ is even and $t \neq r_C$, then $t$ is in a red acyclic component $K$ having $d - zk$ edges and containing a red path $Q_t$ of length $z+1$ ending in $t$, and $K - t$ consists of components with at most $z$ edges.
\end{itemize}
Furthermore, all edges of a red component in $C$ of a vertex $t \in V(T_C)$ are oriented towards $t$ in the example colouring. Note that the paths $Q_t$ and the diameter of $K$ in the second and third case are well-defined by Lemma \ref{lemma:dLargerThanLKPlus1}.

For vertices $t \in V(T_C)$, we again denote by $P_t$ the red component containing $t$ in the example colouring. Note that when creating $C$, there are many possible configurations for the red component $P_t$ if $depth(t)$ is even.
We obtain $G$ by taking $p \geq kD + k^2 + 1$ disjoint copies of $C$ and adding a set $S = \{s_{1},\ldots,s_{k}\}$ of new vertices and
for every colour $i \in \{1, \ldots, k\}$ and every vertex $x$ of a colourful tree we add an edge $xs_i$ \textemdash{} oriented towards $s_i$ in the example colouring \textemdash{} if $x$ does not already have an outgoing edge coloured $i$.

\subsection{Bounding the Diameter}
In this subsection, we prove the theorem below.
\begin{thm} \label{thm:z:lowerBoundDiam}
    In any decomposition of $G$ into $k$ blue pseudoforests and one red pseudoforest where every vertex  has at most $D$ incident red edges, there is a red component $K$ with $e(K) \geq d - z(k-1) + 1$ that has diameter at least $2(z+1)$.
\end{thm}
For the rest of the subsection we assume that there is a colouring $f$ 
of $G$ contradicting Theorem \ref{thm:z:lowerBoundDiam}. 
The following lemma and corollary can be proven analogously to Lemma \ref{lemma:colourfulEdgesToS} and Corollary \ref{cor:forcedColourOfEdges}.

\begin{lemma} \label{lemma:z:colourfulEdgesToS}
    There is a colourful tree $C$ such that every edge of $E(C, S)$ is coloured blue in $f$. Furthermore, there is no monochromatic $S$-$C$-$S$-path.
\end{lemma}

For the rest of the subsection, we let $C$ be a colourful tree satisfying Lemma \ref{lemma:z:colourfulEdgesToS}.

\begin{cor} \label{cor:z:forcedColourOfEdges}
    In $f$, any vertex of $C$ that is not an inner vertex of $T_C$ has a $b$-coloured edge to $S$ for every $b \in \{1, \ldots, k\}$ and any vertex of $T_C$ with odd depth has $k-1$ incident edges to $S$ in pairwise distinct colours of $\{1, \ldots, k\}$. 
    Every edge of $E(C)$ that is not incident to an inner vertex of $T_C$ is coloured red in $f$. Furthermore, every red component containing at least one vertex of $C$ is acyclic in $f$.
\end{cor}

For each vertex $t \in V(T_C)$, let $K_t$ be the subgraph induced by  the vertices $v$ in the same red component in $f$ as $t$ and with $depth_C(v) \geq depth_C(t)$.  Note that $K_t$ is connected and all of its edges are coloured red in $f$. The following lemma bounds $e(K_t)$, and will be used to bound $e(K_{r_C})$ later. Recall that by Corollary \ref{cor:z:forcedColourOfEdges}, every edge not incident to an inner vertex of $T_C$ is red. The main idea below is that if an edge from $t$ to $V(C) \setminus V(T_C)$ is coloured blue in $f$ (thus not matching the example colouring), then we will argue a corresponding edge from $t$ to one of its children $t'$ in $T_C$ is coloured red in $f$; which, together with showing $t'$ is the end of a low $z$-path, will be enough to bound $e(K_t)$. Note that the diameter bounds proceed similarly to the previous section, and the new argument is in bounding the number of edges.

\begin{lemma} \label{lemma:z:inductionOverNonRootVertices}
    If $t \in V(T_C)$ with $t \neq r_C$, then $t$ has a $b$-coloured $S$-path for any $b \in \{1, \ldots, k\}$ and $t$ is the end of a low $z$-path $P$. If $depth(t)$ is even, then $e(P) \geq z+1$ and $e(K_t) \geq d - zk$.
\end{lemma}
\begin{proof}
    The lemma is clear for any leaf of $T_C$ due to Corollary \ref{cor:z:forcedColourOfEdges}. Now, let $t$ have even depth and suppose that the lemma is true for the children of $t$ in $T_C$. 
    All of the children of $t$ in $C$ have a $b$-coloured $S$-path for any $b \in \{1, \ldots, k\}$. Thus, at most $k$ edges from $t$ to the children of $t$ in $C$ can be blue or otherwise there is a monochromatic $S$-$C$-$S$-path, contradicting Lemma \ref{lemma:z:colourfulEdgesToS}.  At least $k+1$ children of $t$ in $C$ are ends of low $z$-paths ($k$ of them are children of $t$ in $T_C$ and the other one is on $Q_t$; recall by Corollary \ref{cor:z:forcedColourOfEdges} the edges of $Q_t$ not incident with $t$ are red). Let us call the set of edges to these $k+1$ children $E_t$.  Since at least one of the edges of $E_t$ is red, $t$ is the end of a low $(z+1)$-path in $f$.
    
    Recall that each of the components of $P_t - t$ has at most $z$ edges, the components of the $k+1$ children of $t$ which are not in $P_t$ have at least $z$ edges, and since at most $k$ edges to the children of $t$ in $C$ are blue (like in the example colouring), we have $e(K_t) \geq e(P_t) = d - zk$.
    
    If fewer than $k$ edges of $E_t$ were blue in $f$, then we had $diam(K_t) \geq 2(z+1)$ and $e(K_t) \geq e(P_t) + (z+1) \geq d - z(k-1) + 1$, which is a contradiction. Thus, $t$ has a $b$-coloured $S$-path for each $b \in \{1, \ldots, k\}$.

    Now, let $t$ have odd depth and not be a leaf and again suppose that the lemma is true for the only child $t'$ of $t$ in $T_C$.  Let $u$ be the other child of $t$ in $C$, which is the neighbour of $t$ in $P_t$. Note that $t'$ has a low $(z+1)$-path and $u$ has a low $(z-1)$-path.
    Since $t$ has $k-1$ $S$-paths of length $1$, at most one of the edges to the children of $t$ can be blue. 
    If both edges were red, then we would have $e(K_t) \geq (d - zk) + (z-1) + 2 = d - z(k-1) + 1$ and $diam(K_t) \geq (z+1) + 2 + (z-1) =  2(z+1)$, which is contrary to our assumptions. The lemma follows.
\end{proof}

\begin{lemma} \label{lemma:z:largeDiamOnRoot}
    We have $e(K_{r_C}) \geq d - z(k-1) + 1$ and $diam(K_{r_C}) \geq 2(z+1)$, which is contradictory to our assumptions.
\end{lemma}
\begin{proof}
    The proof is analogous to the case in the proof of Lemma \ref{lemma:z:inductionOverNonRootVertices} where $t$ has even depth, but this time $|E_t| = k+2$, which forces the statement of the lemma.
\end{proof}

\subsection{Bounding the Fractional Arboricity}

\begin{thm} \label{thm:z:lowerBoundFracArb}
    The fractional arboricity of $G$ is less than $k + \densNDT + \epsilon$. Thus the maximum average degree of $G$ is less than $2\left(k + \densNDT + \epsilon \right)$.
\end{thm}

We assume that there is a subgraph $H$ with $\frac{e(H)}{v(H) - 1} \geq k + \densNDT + \epsilon$ for the rest of the subsection. Let $\mathcal X$ be the set of all $X \subseteq V(G)$ having $S \subsetneq X$ and 
\[\frac{e(X)}{|X \setminus S|} > k + \beta,\]
where
\[
    \beta = \frac{d + \frac{2(z+1)}{\delta + 1}}{d + k + 1 + \frac{2(z+1)}{\delta + 1}}.
\]

\begin{lemma}
    $\frac{e(G)}{|V(G) \setminus S|} = k + \beta$ and thus $V(G) \not\in \mathcal X$.
\end{lemma}
\begin{proof}
    In any colourful tree $C$ there are $\frac{\delta - 1}{2}$ red components $P_t$ in the example colouring, where $depth_{T_C}(t)$ is even and non-zero, and there are $k \frac{\delta + 1}{2}$ such components for which $depth_{T_C}(t)$ is odd. Moreover, every vertex in a colourful tree has exactly one outgoing $i$-edge for each $i \in \{1,2, \dots, k\}$ in the example colouring.
    We have
    \begin{align*}
        \frac{e(G)}{|V(G) \setminus S|}
        &= k + \frac{
            \numSimpleComponents \big(
                d - z(k-1) + 1 + \frac{\delta - 1}{2} (d - zk) + k \frac{\delta + 1}{2} z
            \big)
        }{
            \numSimpleComponents \big(
                d - z(k-1) + 2 + \frac{\delta - 1}{2} (d - zk + 1) + k \frac{\delta + 1}{2} (z + 1)
            \big)
        } \\
        &= k + \frac{\frac{\delta + 1}{2} d + (z+1)}{\frac{\delta + 1}{2} (d+k+1) + (z+1)} \\
        &= k + \beta.
    \end{align*}
\end{proof}

Immediately from the definitions it follows that:

\begin{obs}
     \[
        \densNDT
        < \beta 
        < \densNDT + \epsilon
     \]
\end{obs}

\begin{cor}
    If Theorem \ref{thm:z:lowerBoundFracArb} is false, then $\mathcal X$ is not empty.
\end{cor}
\begin{proof}
    Analogously to Corollary \ref{cor:connectionWeirdDensityToFracArb}.
\end{proof}

Let $C$ be a colourful tree and $t \in V(T_C)$ for the rest of the subsection and let $\ell := \floor{\frac{d-1}{k+1}}$. Note that since $z \leq \ell$, we have $\frac{z}{z+1} \leq \densNDT$.

\begin{lemma}   \label{lemma:z:wholePInX}
    If $X \in \mathcal X$ and $t \in X$, then $X \cup V(P_t) \in \mathcal X$.
\end{lemma}
\begin{proof}
    Analogously to Lemma \ref{lemma:wholePInX}.
\end{proof}

\begin{lemma}
    If $X \in \mathcal X$ and $t \not\in X$, then $X \setminus V(P_t) \in \mathcal X$.
\end{lemma}
\begin{proof}
    Let $P' = X \cap V(P_t) \neq \varnothing$. When subtracting $P'$ from $X$ the induced subgraph will lose $k|P'|$ edges of $E(P', S)$ and another $e(P')$ edges. Let $c$ be the number of components of $G[P']$. By construction every component of $G[P']$ has at most $z$ edges. 
    We have
    \[
        \frac{k|P'| + e(P')}{|P'|}
        \leq k + \frac{cz}{c(z+1)} 
        < k + \beta.
    \]
\end{proof}

Let $V_i$, where $i \in \{z, d - zk, d - z(k-1) + 1\}$, be the set of vertices $t \in V(T_C)$ such that $e(P_t) = i$. Note that these three sets are pairwise disjoint, since $z < d - zk$.

\begin{lemma} \label{lemma:z:ellHasDegkPlusTwo}
    If $X \in \mathcal X$, $t \in X \cap V_z$ and at most  one of the two neighbours of $t$ in $T_C$ are in $X$, then $X \setminus V(P_t) \in \mathcal X$.
\end{lemma}
\begin{proof}
    Analogously to Lemma \ref{lemma:ellHasDegkPlusTwo} using $z \leq \ell$.
\end{proof}

By repeated application of the last three lemmas we get the following corollary:

\begin{cor}
    If $\mathcal X \neq \varnothing$, then there is an $X \in \mathcal X$ such that for every $t \in V(T_C) \cap X$ we have $V(P_t) \subseteq X$ and if additionally $t \in V_z$, then $deg_{X}(t) = k+2$.

\end{cor}

Like in Section 3 we will prove this corollary for all the vertices of $T_C$ now.

\begin{lemma} \label{lemma:z:noPerfectKAryTrees}
    Let $X \in \mathcal X$ such that 
    for every $t \in V(T_C) \cap X$ we have $V(P_t) \subseteq X$ and if additionally $t \in V_z$, then $deg_{X}(t) = k+2$.
    Furthermore, suppose there is a vertex $r' \in V(T_C)$ such that all vertices of the tree $T'$ containing the vertices $t \in V(T_C)$ with $depth(t) \geq depth(r')$ are also in $X$, but the parent of $r'$ is not in $X$.
    Let $T'' = T' \cup \bigcup_{t \in V(T')} P_{t}$.
    Then $X \setminus V(T'')\in \mathcal X$.
\end{lemma}
\begin{proof}
    It must be that $r' \in V_{d - zk}$ because of the degree property of $X$.
    Let $T'$ contain $x$ vertices of $V_{d - zk}$ and thus $kx$ vertices of $V_{z}$. 
    We have
    \begin{align*}
        \frac{ e(T'') + e(X \setminus V(T''), V(T''))} { v(T'') }
        &= k + \frac{ x(d - zk) + xkz }{ x(d+1 - zk) + xk(z+1) } \\
        &= k + \densNDT \\
        &< k + \beta.
    \end{align*}
    Thus, $X \setminus V(T'') \in \mathcal X$.
\end{proof}

\begin{lemma} \label{lemma:z:ellPlusOneHasKChildren}
    If $\mathcal X \neq \varnothing$, then there is an $X \in \mathcal X$ where every $t \in V_{z} \cap X$ has $deg_X(t) = k+2$ and for every $t \in X \cap (V_{d-zk} \cup V_{d - z(k-1) + 1})$ the child of $t$ in $T_C$ is also in $X$.
\end{lemma}
\begin{proof}
    We choose $X, t', T'$ and $X'$ like we did in the proof of Lemma \ref{lemma:ellPlusOneHasKChildren}. Let $x := |V_{d-zk} \cap V(T')|$ and thus $|V_{z} \cap V(T')| = kx + |\{t'\}|$. Since $T' \neq T_C$, we have $x < \frac{\delta - 1}{2}$ and thus
    \[
        \frac{e(X')}{|V(X') \setminus S|}
        =
        k + \frac
            {
                (z + |\{t t'\} | ) + x(d - zk) + kxz
            }{
                z+1 + x(d - zk + 1) + kx(z+1)
            }
        = k + \frac
            {
                xd + (z+1)
            }{
                x(d+k+1) + (z+1)
            }\\
        > \beta.
    \]
\end{proof}

\begin{cor}
    $\mathcal X$ is empty. Thus, Theorem \ref{thm:z:lowerBoundFracArb} is true.
\end{cor}

\begin{ack}
We thank the referees for helpful comments which improved the presentation of the paper. 
\end{ack}

\bibliographystyle{plain}
\bibliography{bib-refs.bib}

\begin{thebibliography}{10}

\bibitem{approxArb}
Markus Blumenstock and Frank Fischer.
\newblock A constructive arboricity approximation scheme.
\newblock In {\em SOFSEM 2020: Theory and Practice of Computer Science}, pages
  51--63. Springer International Publishing, 2020.

\bibitem{dischargingtutorial}
Daniel~W. Cranston and Douglas~B. West.
\newblock An introduction to the discharging method via graph coloring.
\newblock {\em Discrete Mathematics}, 340(4):766--793, 2017.

\bibitem{matroidndt}
Genghua Fan, Hongbi Jiang, Ping Li, Douglas~B. West, Daqing Yang, and Xuding
  Zhu.
\newblock Extensions of matroid covering and packing.
\newblock {\em European Journal of Combinatorics}, 76:117--122, 2019.

\bibitem{ndttPsfs}
Genghua Fan, Yan Li, Ning Song, and Daqing Yang.
\newblock Decomposing a graph into pseudoforests with one having bounded
  degree.
\newblock {\em Journal of Combinatorial Theory, Series B}, 115:72--95, 2015.

\bibitem{digraphndt}
Hui Gao and Daqing Yang.
\newblock Digraph analogues for the nine dragon tree conjecture.
\newblock {\em Journal of Graph Theory}, 102(3):521--534, 2023.

\bibitem{sndtcPsfs}
Logan Grout and Benjamin Moore.
\newblock The pseudoforest analogue for the strong nine dragon tree conjecture
  is true.
\newblock {\em Journal of Combinatorial Theory, Series B}, 145:433--449, 2020.

\bibitem{hakimi}
Seifollah~Louis Hakimi.
\newblock On the degrees of the vertices of a directed graph.
\newblock {\em Journal of the Franklin Institute}, 279(4):290--308, 1965.

\bibitem{ndtt}
Hongbi Jiang and Daqing Yang.
\newblock Decomposing a graph into forests: The nine dragon tree conjecture is
  true.
\newblock {\em Combinatorica}, pages 1125--1137, 2017.

\bibitem{Kostochkaetal}
Seog-Jin Kim, Alexandr~V. Kostochka, Douglas~B. West, Hehui Wu, and Xuding Zhu.
\newblock Decomposition of sparse graphs into forests and a graph with bounded
  degree.
\newblock {\em Journal of Graph Theory}, 74(4):369--391, 2013.

\bibitem{boundeddiameter}
Martin Merker and Luke Postle.
\newblock Bounded diameter arboricity.
\newblock {\em Journal of Graph Theory}, 90(4):629--641, 2019.

\bibitem{miesMoore}
Sebastian Mies and Benjamin Moore.
\newblock The strong nine dragon tree conjecture is true for $d \leq k+1$.
\newblock {\em Combinatorica}, 2023.
\newblock Appeared online.

\bibitem{ndtConjs}
Mickael Montassier, Patrice {Ossona de Mendez}, André Raspaud, and Xuding Zhu.
\newblock Decomposing a graph into forests.
\newblock {\em Journal of Combinatorial Theory, Series B}, 102(1):38--52, 2012.

\bibitem{mousaviHaji}
Seyyed~Ramin Mousavi~Haji.
\newblock Thin trees in some families of graphs.
\newblock Master's thesis, University of Waterloo, 2018.

\bibitem{nash}
Crispin St. J.~A. Nash-Williams.
\newblock Decomposition of finite graphs into forests.
\newblock {\em Journal of the London Mathematical Society 39.1}, page~12, 1964.

\bibitem{Yangmatching}
Daqing Yang.
\newblock Decomposing a graph into forests and a matching.
\newblock {\em Journal of Combinatorial Theory, Series B}, 131:40 -- 54, 2018.

\bibitem{ZHUgamecol}
Xuding Zhu.
\newblock Refined activation strategy for the marking game.
\newblock {\em Journal of Combinatorial Theory, Series B}, 98(1):1--18, 2008.

\end{thebibliography}
    
\end{document}